\documentclass{amsart}

\usepackage{amsmath,amsthm,amssymb,amsfonts,enumerate}

\numberwithin{equation}{section} 
\theoremstyle{plain}
\newtheorem{theo+}           {Theorem}      [section]
\newtheorem{prop+}  [theo+]  {Proposition}
\newtheorem{coro+}  [theo+]  {Corollary}
\newtheorem{lemm+}  [theo+]  {Lemma}
\newtheorem{defi+}  [theo+]  {Definition}

\theoremstyle{definition}
\newtheorem{exam+}  [theo+]  {Example}
\newtheorem{rema+}  [theo+]  {Remark}

\newenvironment{theorem}{\begin{theo+}}{\end{theo+}}
\newenvironment{proposition}{\begin{prop+}}{\end{prop+}}
\newenvironment{corollary}{\begin{coro+}}{\end{coro+}}
\newenvironment{lemma}{\begin{lemm+}}{\end{lemm+}}
\newenvironment{example}{\begin{exam+}}{\end{exam+}}
\newenvironment{remark}{\begin{rema+}}{\end{rema+}}
\newenvironment{definition}{\begin{defi+}}{\end{defi+}}  

\newcommand{\al}{\alpha}
\newcommand{\be}{\beta}       
\newcommand{\ga}{\gamma}  

\newcommand{\de}{\delta}
\newcommand{\De}{\Delta}
\newcommand{\ep}{\varepsilon}
\newcommand{\la}{\lambda}
\newcommand{\ro}{\rho}
\newcommand{\om}{\omega}

\newcommand{\si}{\sigma}

\newcommand{\La}{\Lambda}
\newcommand{\ka}{\kappa}
\newcommand{\gh}{\mathfrak h}
\newcommand{\gha}{\mathfrak h^\ast}
\newcommand{\mg}{M_{\gh^\ast}}
\newcommand{\mc}{M_{\mathbb C}}
\newcommand{\D}{D_{\gh}}
\newcommand{\I}{I_{\gh}}
\newcommand{\End}{\operatorname{End}}
\newcommand{\id}{\operatorname{id}}
\newcommand{\Hom}{\operatorname{Hom}}

\newcommand{\cop}{\operatorname{cop}}
\newcommand{\Ima}{\operatorname{Im}}
\newcommand{\Ker}{\operatorname{Ker}}
\newcommand{\ot}{\otimes}
\newcommand{\C}{\mathbb C}

\newcommand{\Z}{\mathbb Z}
\newcommand{\Zp}{\mathbb Z_{\geq 0}}
\newcommand{\wh}{\widehat}
\newcommand{\wt}{\widetilde}
\newcommand{\fr}{\mathcal F_R}

\newcommand{\uq}{\mathcal U_q(\mathrm{sl}(2))}
\newcommand{\lx}{\langle}
\newcommand{\rx}{\rangle}

\newcommand{\qb}[2]{\genfrac{[}{]}{0pt}{}{#1}{#2}_{q^2}}

\begin{document}
\baselineskip 18pt
\larger[2]
\title
[Duality for dynamical quantum groups]
{Duality and self-duality\\ for dynamical quantum groups}
\author{Hjalmar Rosengren}
\address
{Department of Mathematics\\ Chalmers University of Technology and G\"oteborg
 University\\SE-412~96 G\"oteborg, Sweden}
\email{hjalmar@math.chalmers.se}

\begin{abstract}
We define a natural  concept of duality for the $\gh$-Hopf algebroids
introduced by Etingof and Varchenko. We prove that the special case
of the trigonometric
 $\mathrm{SL}(2)$ dynamical quantum group is self-dual, and may 
therefore be viewed as a deformation  
both of the function algebra $\mathcal F(\mathrm{SL}(2))$ and of the
enveloping algebra $\mathcal U(\mathrm{sl}(2))$. Matrix elements of
the self-duality in the Peter--Weyl basis are $6j$-symbols; this leads
to a new algebraic interpretation of the hexagon identity or quantum dynamical
Yang--Baxter equation for quantum and classical $6j$-symbols.  
\end{abstract}

\maketitle        

\section{Introduction}  
Quantum groups may be viewed as an algebraic framework for studying
solutions to the  Yang--Baxter equation, which exists in several versions. 
Although the quantum dynamical Yang--Baxter (QDYB) 
equation  is much older (cf.\ \S \ref{ssf} for some historical comments) 
than the perhaps more well-known 
 equation  with spectral parameters, the corresponding quantum groups have only
been studied in recent years. The main reason is perhaps that these objects
are not Hopf algebras, but  more involved structures such as
Hopf algebroids \cite{ev,ev2,x}, weak Hopf algebras \cite{en} and 
quasi-Hopf algebras \cite{bbb,jkos}. 

The dynamical quantum groups of this paper are so called
 ``$\gh$-Hopf algebroids'', a notion
introduced  by Etingof and Varchenko \cite{ev}, motivated by earlier
work of Felder and Varchenko \cite{f,fv}. They are constructed from 
solutions to the QDYB equation in a manner analogous to the 
Faddeev--Reshetikhin--Sklyanin--Takthajan (FRST) construction. 

In \cite{ev2} it was suggested that a large class of
 dynamical quantum groups  are  self-dual. This 
should mean that they are analogues both of function algebras and of
envelopping algebras, a rather intriguing fact. 
However, to quote \cite{en}: ``It is not very convenient to formulate such a
statement precisely, because of difficulties with the notion 
of a dual Hopf algebroid.'' It is the purpose of this paper 
to resolve these difficulties. We  show how to formulate a duality theory 
for $\gh$-Hopf algebroids and
 prove  self-duality in the $\mathrm{SL}(2)$ case.
The authors of \cite{en} chose a different
approach and obtained a  general self-duality theorem
for dynamical quantum groups within the framework of weak Hopf algebras. 
Though it may be possible to transfer the results of \cite{en} to 
$\gh$-Hopf algebroids, the approach of this paper is completely different
and should be of independent interest. 

It is interesting to compare dynamical quantum groups with another type
of self-dual quantum groups: the braided groups of Majid
\cite{m}. For braided groups, the self-duality degenerates 
as the deformation parameter $q\rightarrow 1$. By contrast,
for dynamical
quantum groups there is nothing special (from an algebraic viewpoint)
about the case $q=1$; in particular  self-duality still holds. 

Let us summarize the contents of the paper. In \S 2 we recall the necessary
algebraic background. References for this material are \cite{ev} and 
\cite{kr}. In \S 3 we show that there is a working duality theory
for $\gh$-Hopf  algebroids. This leads naturally to the definition of
a cobraiding on an $\gh$-Hopf  algebroid.
The main result in this section is Corollary \ref{rp}, which shows that 
applying the 
generalized FRST construction  of \cite{ev} to a dynamical $R$-matrix
automatically gives a \emph{cobraided} $\gh$-bialgebroid.
A corresponding statement is true for Hopf algebras, but in interesting 
examples the cobraiding will  be degenerate. In \S 4 we prove that 
for the trigonometric $\mathrm{SL}(2)$
dynamical quantum group, the cobraiding is ``almost'' non-degenerate.
The radical of the cobraiding serves to eliminate those representations
which correspond to some covering group rather than $\mathrm{SL}(2)$.
Finally we show that the matrix elements of the cobraiding
in the Peter--Weyl basis are quantum (classical if $q=1$) $6j$-symbols. 
This allows us to recover the first known instance of the QDYB equation: the
hexagon identity for $6j$-symbols found by Wigner in 1940. 

{\bf Acknowledgment:} I would like to thank
 Erik Koelink for valuable discussions.

\section{Preliminaries}
\subsection{$\gh$-algebra}
\label{has}

Throughout the paper, $\gha$ will be a finite-dimensional complex 
vector space. In the context of dynamical quantum groups, it appears as
the dual of a Cartan subalgebra of the corresponding Lie algebra.
We denote by $\mg$ the field of meromorphic functions on $\gha$.

An \emph{ $\gh$-prealgebra} $A$ is a complex vector space, equipped with a 
decomposition $A=\bigoplus_{\al,\be\in\gha} A_{\al\be}$ 
and two left actions $\mu_l,\,\mu_r:\,\mg\rightarrow\End_\mathbb C(A)$
(the left and right moment maps)
which preserve the bigrading,  such that the images of $\mu_l$ and
$\mu_r$ commute.  A homomorphism of $\gh$-prealgebras is
a linear map which preserves the moment maps and the bigrading.
It is convenient to   introduce two right actions by
\begin{equation}\label{dr}\mu_l(f)a=a\mu_l(T_\al f),
\qquad \mu_r(f)a=a\mu_r(T_\be f), \qquad a\in A_{\al\be},\ f\in \mg,
\end{equation}
where $T_\al$ denotes the automorphism $T_\al f(\la)=f(\la+\al)$ of $\mg$.

We need two different tensor products on
$\gh$-prealgebras. The first one, denoted  $A\widehat\otimes B$,
equals $A\ot_{\mathbb C} B$ modulo
the relations
$$a\mu_l^A(f)\otimes b=a\otimes\mu_l^B(f)b,
\qquad  a\mu_r^A(f)\otimes b=a\otimes\mu_r^B(f)b.$$
The bigrading $A_{\al\be}\wh\otimes
A_{\ga\de}\subseteq(A\widehat\otimes B)_{\al+\ga,\,\be+\de}$
and the moment maps 
$$\mu_l^{A\widehat\otimes B}(f)(a\otimes b)=\mu_l^A(f)a\otimes
b,\qquad
\mu_r^{A\widehat\otimes B}(f)(a\otimes b)=\mu_r^A(f) a\otimes b$$
 make $A\widehat\otimes B$ an
$\gh$-prealgebra.

Another kind of tensor product, denoted
 $A\widetilde\otimes B$, equals $\bigoplus_{\al\be\ga}
 A_{\al\ga}\otimes_{\mathbb C} B_{\ga\be}$ modulo the relations
$$\mu_r^A(f)a\otimes b=a\otimes \mu_l^B(f) b,
\qquad a\in A,\ b\in B,\ f\in
\mg.$$
The bigrading 
$A_{\al\be}\wt\ot B_{\be\ga}\subseteq(A\wt\ot B)_{\al\ga}$
and the moment maps
$$\mu_l^{A\widetilde\otimes B}(f)(a\ot b)=\mu_l^A(f)a\otimes b,
\qquad \mu_r^{A\widetilde\otimes
B}(f)(a\ot b)=a\otimes\mu_r^B(f)b$$
make $A\widetilde\otimes B$  an $\gh$-prealgebra.

An \emph{$\gh$-algebra} is an $\gh$-prealgebra which is also an associative 
algebra with $1$. It is required  that the decomposition is a bigrading:
$A_{\al\be}A_{\ga\de}\subseteq A_{\al+\ga,\be+\de}$. Considering
 $\mu_l$ and $\mu_r$ as algebra embeddings $\mg\rightarrow A_{00}$
through $\mu_l(f)=\mu_l(f)1$,  it is moreover required that
\eqref{dr} hold 
as  relations in the algebra. A homomorphism of $\gh$-algebras is an
$\gh$-prealgebra homomorphism which is also an algebra homomorphism.
If $A$ and $B$ are $\gh$-algebras, then so is $A\widetilde\otimes B$, 
with the multiplication $(a\otimes b)(c\otimes d)=ac\otimes bd$.

We denote by $D_\gh$ the algebra of difference operators on $\mg$,
consisting of finite sums
$\sum_if_i\,T_{\be_i}$, $f_i\in \mg$, $\be_i\in\gh^\ast$.
This is an $\gh$-algebra with the bigrading
defined by $f\,T_{-\be}\in (D_\gh)_{\be\be}$ and both moment maps equal
to  the natural
embedding. For any $\gh$-algebra $A$, there are canonical $\gh$-algebra
isomorphisms 
$A\simeq A\widetilde\otimes D_\gh\simeq D_\gh\widetilde\otimes A$, defined by
\begin{equation}\label{can}x\simeq x\otimes T_{-\be}\simeq
T_{-\al}\otimes x,\qquad x\in A_{\al\be}.\end{equation}

An \emph{$\gh$-coalgebroid} is an $\gh$-prealgebra equipped with  two 
$\gh$-prealgebra
homomorphisms, $\De:\, A\rightarrow A\widetilde\otimes A$ (the coproduct) and
 $\ep:\,A\rightarrow D_{\gh}$ (the counit), such that 
$(\De\otimes\id)\circ\De=(\id\otimes\De)\circ\De$ and, under the 
identifications
(\ref{can}),  $(\ep\otimes\id)\circ\De=(\id\otimes\,\ep)\circ\De=\id$.
 An $\gh$-coalgebroid homomorphism
$\phi:\,A\rightarrow B$ is an $\gh$-prealgebra homomorphism with
$(\phi\otimes\phi)\circ\De^A=\De^B\circ\phi$, $\ep^B\circ\phi=\ep^A$.

If $A$ and $B$ are $\gh$-coalgebroids, then 
$\De^{A\wh\ot B}=\si_{23}\circ(\De^A\ot \De^B)$ and
$\ep^{A\widehat\otimes B}(a\otimes b)=\ep^A(a)\circ\ep^B(b)$
 define an $\gh$-coalgebroid structure on
$A\widehat\otimes B$. Here, $\si_{23}$ is a quotient of the map
$a\ot b\ot c\ot d\mapsto a\ot c\ot b\ot d$.

For completeness we note here (this is not in the references)
that the product $\wh\ot$ has a unit element. Namely, let
 $\I=\mg\ot_{\C}\mg$ with  the bigrading $\I=(\I)_{00}$, the moment maps
$\mu_l(f)=f\ot 1$, $\mu_r(f)=1\ot f$, the coproduct
$\De(f\ot g)= f\ot 1\ot 1\ot g$ and the counit $\ep(f\ot g)=fg$. 
Then $\I$ is an $\gh$-coalgebroid. 
Moreover, for any $\gh$-coalgebroid $A$ there are canonical
$\gh$-coalgebroid isomorphisms $A\simeq A\wh\ot \I\simeq \I\wh\ot A$,
defined by
\begin{equation}\label{can2}
x\simeq x\ot 1\ot 1\simeq 1\ot 1\ot x,\qquad x\in A.\end{equation}

An \emph{$\gh$-bialgebroid} is an $\gh$-algebra which is also an 
$\gh$-coalgebroid, such that $\De$ and $\ep$ are $\gh$-algebra homomorphisms.
It follows that the multiplication in $A$ factors to an
$\gh$-coalgebroid homomorphism $A\wh\ot A\rightarrow A$ and that
$f\ot g\mapsto\mu_l(f)\mu_r(g)$ defines an $\gh$-coalgebroid
homomorphism $\I\rightarrow A$.

In particular, if $A$ and $B$ are $\gh$-bialgebroids,
then $A\wh\ot B$ is an $\gh$-coalgebroid  and
$A\wt\ot B$ an $\gh$-algebra. We stress that
none of these spaces carry a natural $\gh$-bialgebroid structure.

An \emph{$\gh$-Hopf algebroid} is an $\gh$-bialgebroid $A$ equipped
with a map $S\in \End_\C(A)$ (the antipode), such that
\begin{gather}\label{sf}S(\mu_r(f)a)=S(a)\mu_l(f),\qquad 
S(a\mu_l(f))=\mu_r(f)S(a),
\qquad a\in A,\ f\in \mg,\\
\begin{split}m\circ(\id\otimes\, S)\circ\De(a)&
=\mu_l(\ep(a)1),\qquad a\in A,\\
\label{anti}m\circ(S\otimes \id)\circ\De(a)&=\mu_r(T_\al(\ep(a)1)),
\qquad a\in A_{\al\be},
\end{split}\end{gather}
where $m$ denotes multiplication and where $\ep(a)1$ is the action of the 
difference operator $\ep(a)$ on the function $1\in\mg$. 
In \cite{kr} it was proved that 
$S$ is unique and is an antihomomorphism (in a natural sense)
of $\gh$-bialgebroids.

The $\gh$-algebra $D_{\gh}$ is an $\gh$-Hopf algebroid with $\De^{D_\gh}$ the
canonical isomorphism as in \eqref{can}, 
$\ep^{D_{\gh}}$ the identity map and the 
antipode defined by $S^{D_{\gh}}(f T_\al)=(T_{-\al}f)T_{-\al}$, or equivalently
$S^{D_{\gh}}(x)=T_{\al}\circ x\circ T_{\al}$, 
$x\in (D_\gh)_{\al\al}$.  

The $\gh$-coalgebroid  $\I$ is an $\gh$-Hopf algebroid with the product 
$(f\ot g)(h\ot k)=fh\ot gk$ and the
antipode $S(f\ot g)=g\ot f$.

\subsection{Representation theory}

When considering representations and corepresentations 
of $\gh$-algebroids, the representation
spaces must also carry an $\gh$-structure. More precisely, 
by an \emph{$\gh$-space} we mean an $\gh^\ast$-graded vector space over $\mg$,
$V=\bigoplus_{\al\in\gha} V_\al$. 
A morphism of $\gh$-spaces is a
grade-preserving $\mg$-linear map. In analogy with \eqref{dr}, we will  write
$$vf=T_{-\al}(f)v,\qquad f\in\mg,\ v\in V_\al.$$

  If $A$ is an $\gh$-prealgebra
 and $V$ an $\gh$-space, we define 
 $A\widetilde\otimes V=\bigoplus_{\al\be} A_{\al\be}\otimes_{\C}V_\be$
    modulo the relations
 $$\mu_r^A(f)a\otimes v=a\otimes fv.$$ 
 The grading
  $A_{\al\be}\wt\otimes V_\be\subseteq(A\widetilde\otimes V)_\al$ and the
$\mg$-linear structure 
 $f(a\otimes v)=\mu_l^A(f)a\otimes v$
  make $A\widetilde\otimes V$ into
 an $\gh$-space. A \emph{corepresentation}
 of an $\gh$-coalgebroid $A$
 on an  $\gh$-space
$V$ is an $\gh$-space morphism $\pi:\,V\rightarrow A\widetilde\otimes V$ 
such that
$(\De\otimes\id)\circ\pi=(\id\otimes\,\pi)\circ\pi$,
$(\ep\otimes\id)\circ\pi=\id.$
The second equality is in the sense of the  isomorphism 
$D_{\gh}\widetilde\otimes V\simeq V$ defined
by $f\,T_{-\al}\otimes v\simeq fv$, $f\in\mg$, $v\in V_\al$.
If $\{v_k\}_k$ is a homogeneous basis (over $\mg$) of $V$, 
$v_k\in V_{\om(k)}$, the matrix elements $t_{kj}\in A_{\om(k),\,\om(j)}$ of 
$\pi$ with respect to this basis are given by
\begin{equation}\label{cme}\pi(v_k)=\sum_j t_{kj}\ot v_j.\end{equation} 
They satisfy
\begin{equation}\label{dt}\De(t_{kj})=\sum_l t_{kl}\ot t_{lj},
\qquad \ep(t_{kj})=\de_{kj}\,T_{-\om(k)}.\end{equation} 
An intertwiner $\Phi:\,V\rightarrow W$ of corepresentations is
an $\gh$-space morphism which satisfies 
$\pi_W\circ\Phi=(\id\ot\Phi)\circ \pi_V$. Given two bases $\{v_k\}_k$, 
$\{w_k\}_k$ of $V$ and $W$, respectively, and writing
 $\Phi(v_k)=\sum_j\Phi_{kj} w_j$, the intertwining property may be written
\begin{equation}\label{ix}
\sum_j\mu_r(\Phi_{jl})t_{kj}^V=\sum_j\mu_l(\Phi_{kj})t_{jl}^W
\qquad \text{for all } k \text{ and } l. \end{equation}

Given two $\gh$-spaces $V$ and $W$, their tensor product
$V\wh\ot W$ is defined as $V\ot W$ modulo the relations 
$$vf\ot w=v\ot fw$$
and equipped with the $\gh$-space structure $V_\al\wh\ot W_\be\subseteq
(V\wh\ot W)_{\al+\be}$, $f(v\ot w)=fv\ot w$. 
If $V$ and $W$ are 
corepresentations of an $\gh$-bialgebroid, then so is $V\wh\ot W$. 
Explicitly, if $t_{kj}^V$ and $t_{kj}^W$ are
matrix elements as in \eqref{ix}, then
\begin{equation}\label{tpx}
\pi_{V\wh\ot W}(v_j\ot w_k)=\sum_{lm}t_{jl}^V\,t_{km}^W\ot v_l\ot w_m.
\end{equation}

Next we turn to representations. For $V$ an $\gh$-space,  let 
 $(D_{V})_{\al\be}$ be the space of $\mathbb C$-linear operators $U$
on  $V$
such that $U(gv)=T_{-\be}(g)U(v)$ and $U(V_\ga)\subseteq V_{\ga+\be-\al}$ for 
all $g\in\mg$, $v\in V$ and $\ga\in\gh^\ast$. Then the space
  $D_{V}=\bigoplus_{\al\be\in\gh^\ast}
(D_{V})_{\al\be}$ is an $\gh$-algebra with
 moment maps 
\begin{equation}\label{md}\mu_l(f)(v)=vf,\qquad \mu_r(f)(v)=f v.\end{equation}
An  \emph{$\gh$-representation} 
 of an $\gh$-algebra $A$ on  $V$ is
an $\gh$-algebra homomorphism $A\rightarrow D_{V}$. (In \cite{ev} this
was called a dynamical representation. The  attribute
$\gh$- or dynamical should remind us that we are not merely representing the 
 underlying associative algebra, but also
the $\gh$-structure.)    
An intertwiner  $\Phi:\,V\rightarrow W$ of $\gh$-representations is
an  $\gh$-space morphism which satisfies $\Phi\circ\pi_V(a)=\pi_W(a)\circ\Phi$
for all $a\in A$.

 Given two $\gh$-representations $V$ and $W$ of an $\gh$-bialgebroid $A$,
$V\wh\ot W$ is an $\gh$-representation with
$$\pi_{V\wh\ot W}(a)(v\ot w)=\sum_i \pi_V(a_i'')v\ot \pi_W(a_i')w,\qquad
\De(a)=\sum_i a_i'\ot a_i''.$$
Since this is an analogue of the representation $W\ot V$ for bialgebras,
it is perhaps more natural to work with the tensor product opposite 
to $\wh\ot$ (denoted $\bar\ot$ in \cite{ev}), but we avoid that in
the present paper.

\subsection{Generalized FRST construction}
\label{ssf}

Let $X$ be a finite index set, $\om:\,X\rightarrow\gha$ an arbitrary 
function and $R=(R_{xy}^{ab})_{x,y,a,b\in X}$ a matrix with
entries in $\mg$ such that $R_{xy}^{ab}=0$ if $\om(x)+\om(y)\neq\om(a)+\om(b)$.
To this data is associated an $\gh$-bialgebroid $A_R$ generated by 
$\{L_{xy}\}_{x,y\in X}$ together with two copies of $\mg$, embedded as
subalgebras. We write the elements of these copies as $f(\la)$, $f(\mu)$,
respectively. The defining relations of $A_R$
are
\begin{equation}\begin{gathered}\label{zr}f(\la)L_{xy}=L_{xy}f(\la+\om(x)),
\qquad
 f(\mu)L_{xy}=L_{xy}f(\mu+\om(y)),\\
 f(\la)g(\mu)=g(\mu)f(\la)\end{gathered}\end{equation}
for $f$, $g\in\mg$, together with 
\begin{equation}\label{rll}\sum_{xy}R_{ac}^{xy}(\la)L_{xb}L_{yd}=
\sum_{xy}R_{xy}^{bd}(\mu)L_{cy}L_{ax},\qquad a,b,c,d\in X.\end{equation}
The bigrading on $A_R$ is defined by $L_{xy}\in A_{\om(x),\,\om(y)}$, 
$f(\la),\,
f(\mu)\in A_{00}$, and the moment maps by $\mu_l(f)=f(\la)$,
$\mu_r(f)=f(\mu)$.  The coproduct and  counit are defined by
\begin{gather*}\De(L_{ab})=\sum_{x\in X}L_{ax}\otimes L_{xb},
\qquad \De(f(\la))=f(\la)\otimes 1,\qquad \De(f(\mu))=1\otimes f(\mu),\\
\ep(L_{ab})=\de_{ab}\,T_{-\om(a)},
\qquad \ep(f(\la))=\ep(f(\mu))=f.\end{gather*}

Let $V$ be a complex vector space with basis $\{v_x\}_{x\in X}$, viewed as
an $\gh$-module through $v_x\in V_{\om(x)}$. We  identify 
$R$ with
 the  meromorphic function $\gh^\ast\rightarrow\End_\gh(V\ot V)$
defined by 
\begin{equation}\label{me}
R(\la)(v_a\otimes v_b)=\sum_{xy}R_{xy}^{ab}(\la)\,v_x\otimes v_y.\end{equation}
Then $R$ is called a \emph{dynamical $R$-matrix} if it satisfies the quantum
dynamical Yang--Baxter (QDYB) equation 
\begin{equation}\label{dy}
R^{12}(\la-h^{(3)})R^{13}(\la)R^{23}(\la-h^{(1)})=R^{23}(\la)
R^{13}(\la-h^{(2)})R^{12}(\la).\end{equation}
This is an identity in the algebra of meromorphic functions $\gh^{\ast}
\rightarrow\End(V\otimes V\otimes V)$. Here,
 $h$ indicates the action of $\gh$, and the upper
indices refer to the factors in the tensor product. For instance, 
$R^{12}(\la-h^{(3)})$ denotes the operator
$$R^{12}(\la-h^{(3)})(u\otimes v\otimes w)=(R(\la-\mu)(u\otimes v))\otimes
w,\qquad w\in V_\mu.$$
Evaluating both sides of \eqref{dy} on a tensor product $v_a\ot v_b\ot v_c$
and identifying the coefficient of $v_d\ot v_e\ot v_f$ gives the
expression for the QDYB equation
in terms of the matrix elements of $R$:
\begin{equation}\label{dyx}\sum_{xyz} R^{xy}_{de}(\la-\om(f))
\,R^{az}_{xf}(\la)\,R^{bc}_{yz}(\la-\om(a))
=\sum_{xyz}R^{yz}_{ef}(\la) \,R^{xc}_{dz}(\la-\om(y))\,R^{ab}_{xy}(\la).
\end{equation}
Although the FRST construction works for general $R$,  the case
when $R$ is a dynamical $R$-matrix is the most interesting
(cf.~\cite{ev} and Corollary \ref{rp} below).

The dynamical Yang--Baxter equation first occurred in the work of
Wigner \cite{w}. The matrix elements $R_{xy}^{ab}$ are then $6j$-symbols and
each side of \eqref{dyx} a $9j$-symbol, so the QDYB equation expresses
a symmetry of the $9j$-symbol.
 Although it is common knowledge that Wigner's identity is
a kind of Yang--Baxter equation, it seems that it was first written down
in the form \eqref{dy} by Gervais and Neveu \cite{gn} (where it arose 
independently and in 
a different context). Therefore, the QDYB equation has also been called
the Gervais--Neveu equation. 

\subsection{The $\mathrm{SL}(2)$ dynamical quantum group}

Our main example is the trigonometric 
$\mathrm{SL}(2)$ dynamical quantum group. 
In this example
$\gh=\gh^\ast=\C$, $X=\{+,-\}$, $\om(\pm)=\pm 1$ and $R$
is the dynamical $R$-matrix
\begin{equation}\label{r}\left(\begin{matrix}
 R_{++}^{++} & R_{++}^{+-} & R_{++}^{-+} & R_{++}^{--}\\
 R_{+-}^{++} & R_{+-}^{+-} & R_{+-}^{-+} & R_{+-}^{--}\\
 R_{-+}^{++} & R_{-+}^{+-} & R_{-+}^{-+} & R_{-+}^{--}\\
R_{--}^{++} & R_{--}^{+-} & R_{--}^{-+} & R_{--}^{--}\end{matrix}\right)
=\left(\begin{matrix}q&0&0&0\\
0&1&\frac{q^{-1}-q}{q^{2(\lambda+1)}-1}&0\\
0&\frac{q^{-1}-q}{q^{-2(\lambda+1)}-1}&\frac{(q^{2(\lambda+1)}-q^2)
(q^{2(\lambda+1)}-q^{-2})}{(q^{2(\lambda+1)}-1)^2}&0\\
0&0&0&q
 \end{matrix}\right),
\end{equation}
with $q$ a fixed parameter, $0<q<1$. This is the $R$-matrix arising from
$6j$-symbols of the standard quantum group $\mathcal U_q(\mathrm{sl}(2))$
evaluated in  a two-dimensional representation.
It can also be obtained by a twisting construction from the 
$R$-matrix of $\mathcal U_q(\mathrm{sl}(2))$, using certain 
``dynamical boundaries''  discovered by Babelon \cite{b,bbb}; cf.\ also
\cite{r}.  

We write  the generators of the corresponding $\gh$-bialgebroid 
 as
\begin{equation}\label{lg}\al=L_{++},\qquad \be=L_{+-},\qquad
\ga=L_{-+},\qquad \de=L_{--}.\end{equation} 
In terms of the auxiliary functions
\begin{align*}F(\la)&=\frac{q^{2(\la+1)}-q^{-2}}{q^{2(\la+1)}-1},\qquad
G(\la)=\frac{(q^{2(\la+1)}-q^2)(q^{2(\la+1)}-q^{-2})}{(q^{2(\la+1)}-1)^2},\\
H(\la,\mu)&=\frac{(q-q^{-1})(q^{2(\lambda+\mu+2)}-1)}
{(q^{2(\la+1)}-1)(q^{2(\mu+1)}-1)},\qquad
I(\la,\mu)=\frac{(q-q^{-1})(q^{2(\mu+1)}-q^{2(\la+1)})}
{(q^{2(\la+1)}-1)(q^{2(\mu+1)}-1)},\end{align*}   
the sixteen relations \eqref{rll} reduce to the six independent
equations 
\begin{gather*}
\al\be=qF(\mu-1)\be\al,\qquad \al\ga=qF(\la)\ga\al,\qquad \be\de=qF(\la)\de\be,
\qquad \ga\de=qF(\mu-1)\de\ga,\\
\al\de-\de\al=H(\la,\mu)\ga\be,\qquad \be\ga-G(\la)\ga\be
=I(\la,\mu)\al\de.
\end{gather*}

It is possible to obtain an $\gh$-Hopf-algebroid, which we denote
 $\fr=\mathcal F_R(\mathrm{SL}(2))$, 
from this example by
adjoining the relation
\begin{equation}\label{ddr}\al\de-qF(\la)\ga\be=1\end{equation}
and defining the antipode by
$$
S(\al)=\frac{F(\la)}{F(\mu)}\,\de,
\qquad S(\be)=-\frac{q^{-1}}{F(\mu)}\,\be,
\qquad S(\ga)=-qF(\la)\ga,\qquad S(\de)=\al.$$

The Hopf algebra  $\mathcal F_q(\mathrm{SL}(2))$ is recovered as the
formal limit of $\fr$ when the dynamical variables $\la$, 
$\mu\rightarrow-\infty$. We will refer to this as the
\emph{non-dynamical limit}. Another interesting limit is the rational limit
$q\rightarrow 1$. We stress that 
 all our results survive the rational limit. That is, 
the rational $\mathrm{SL}(2)$ dynamical quantum group is an essentially
self-dual Hopf algebroid, which is constructed from the classical $6j$-symbols
of Racah and Wigner (very natural objects in the representation theory of
$\mathrm{SL}(2)$), and in fact provides an alternative algebraic framework 
for deriving their main properties.

An important difference between the dynamical and non-dynamical case is that
$\fr$ has a non-trivial center. In fact, the element
\begin{equation}\label{xi}\Xi=q^{\la-\mu-1}+q^{\mu-\la+1}
-q^{-(\la+\mu+2)}(1-q^{2(\la+1)})(1-q^{2(\mu+1)})\be\ga
\end{equation}
is central. In \cite{kr} it was observed that $\Xi$ plays the role
of Casimir element of $\fr$. We will obtain a precise version of this 
statement in Proposition \ref{ufp} below.

\section{Duality for $\gh$-bialgebroids}

\subsection{Algebraic duals}\label{ssad}

In this section we show that there is a working duality theory for 
$\gh$-bialgebroids. As in the case of Hopf algebras, it is  more
convenient in practice to work with pairings between two
objects than to work directly with the algebraic dual. 
Consequently, the reader may wish to skip this part of the paper
 and pass directly to the 
definition of a pairing in \S \ref{sps}. 

 Let us first try to motivate our construction.
If $A$ is an $\gh$-coalgebroid, we expect the counit $\ep^A$ to be the
unit in the dual $\gh$-algebra $A'$. Thus, we want $A'$ to be
a subspace of $\Hom_\C(A,D_\gh)$. It will be convenient to write
$$\langle a,\phi\rangle=\phi(a),\qquad a\in A,\ \phi\in\Hom_\C(A,D_\gh).$$
 It follows from the
counit axioms that if $\De(a)=\sum_i a_i'\ot a_i''$, $a_i'\in A_{\al\be_i}$,
$a_i''\in A_{{\be_i}\ga}$, then
\begin{equation}\label{eu}\lx x,\phi\rx=\sum_i \lx a_i',\ep\rx
\,T_{\be_i}\,\lx a_i'',\phi\rx=
\sum_i \lx a_i',\phi\rx\,T_{\be_i}\,\lx a_i'', \ep\rx\end{equation}
for suitable $\phi\in\Hom_\C(A,D_\gh)$. By suitable, we mean that the above 
expression should be independent of the choice of representative in
$A\ot A$ for $\De(a)\in A\wt\ot A$, that is, that
\begin{align*}\sum_i\lx\mu_r(f)a_i', \ep\rx\,T_{\be_i}\,\lx a_i'', \phi\rx
&=\sum_i\lx a_i',\ep\rx\,T_{\be_i}\,\lx\mu_l(f)a_i'', \phi\rx,\\
\sum_i \lx \mu_r(f)a_i',\phi\rx\,T_{\be_i}\,\lx a_i'', \ep\rx&=
 \sum_i\lx a_i',   \phi\rx\,T_{\be_i}\,\lx \mu_l(f)a_i'', \ep\rx\end{align*}
for all $f\in\mg$.
A sufficient condition for this to hold is that 
\begin{equation}\label{ap0}
\lx\mu_l(f)a,  \phi\rx=f\circ\lx a,\phi\rx,\qquad 
\lx a\mu_r(f), \phi\rx=\lx a,\phi\rx\circ f,\qquad a\in A,\ f\in\mg. 
\end{equation}
Here and below we write $f$ for the operator  $fT_0\in D_\gh$.

Let us write $A'_{\al\be}$ for the subspace of 
$\Hom_\C(A,D_\gh)$ consisting of elements $\phi$ satisfying \eqref{ap0}
such that 
$$\langle A_{\ga\de},A'_{\al\be}\rangle\subseteq(D_\gh)_{\ga+\be,\de+\al}$$
(in particular, $\langle A_{\ga\de},A'_{\al\be}\rangle=0$ if 
$\ga+\be\neq \de+\al$),
and define
$A'=\bigoplus_{\al\be\in \gh^\ast}A'_{\al\be}$. 
It is  easy to check that $A'$ is closed under the multiplication
\begin{equation}\label{dp}\lx a,\phi\psi\rx=\sum_i \lx a_i', \phi\rx\,
T_{\be_i}\,\lx a_i'', \psi\rx.\end{equation}
It follows from the coproduct axiom that this product is associative, and
from \eqref{eu} that $\ep$ is a unit element. Finally, we introduce the
moment maps
$$\lx a,\mu_l(f)\rx=f\circ \lx a,\ep\rx,
\qquad \lx a,\mu_r(f)\rx=\lx a, \ep\rx\circ f,$$
so that
$$\lx a, \mu_l(f)\phi\rx=f\circ\lx a,\phi\rx=\lx \mu_l(f)a,\phi\rx
,\qquad\lx a,\phi\mu_r(f)\rx=\lx a,\phi\rx\circ f=\lx a\mu_r(f),\phi\rx.$$
It is then easy to check that $A'$ is an $\gh$-algebra. 
We call it the  \emph{dual $\gh$-algebra} of the $\gh$-coalgebroid $A$.

\begin{remark} A  product related to \eqref{dp},
$$\lx a,\phi\star\psi\rx=\sum_i \lx  a_i',\phi\rx\lx a_i'',\psi\rx,$$
which makes sense on a certain subspace of $\Hom_\C(A,B)$ for $B$ an arbitrary
$\gh$-algebra, was used and studied in \cite{kr}. This product has
 a left unit coming from $\mu_r^B$
and a right unit coming from $\mu_l^B$. It is the fact that
$\mu_l^{D_\gh}=\mu_r^{D_\gh}$ that allows us to modify it into a product
 with a two-sided unit. 
\end{remark}

Next we want to define the dual $\gh$-bialgebroid of an $\gh$-bialgebroid.
We will need the following  lemmas. We omit the straight-forward
proof of the first one.

\begin{lemma}\label{dhl}
Let $A$ and $B$ be $\gh$-coalgebroids and $\chi:\,A\rightarrow B$ an
$\gh$-coalgebroid homo\-morphism. Then $\chi'(\phi)=\phi\circ\chi$
defines an $\gh$-algebra homomorphism $\chi':\,B'\rightarrow A'$.
\end{lemma}

\begin{lemma}\label{ld}
Let $A$ be an $\gh$-coalgebroid and let
 $\psi_1,\dots,\psi_n\in A'$ be linearly independent over
$\mu_l^{A'}(\mg)$.
Then there exists $b\in A$ with $\lx b,\psi_i\rx 1=\de_{i1}$.
\end{lemma}

\begin{proof}
We view $A$ and $A'$ as vector spaces over $\mg$ through $\mu_l^A$, 
$\mu_l^{A'}$. Let $X=\bigoplus_{i=1}^n\mg\psi_i$.
Then $\ro(a)(\psi)=\lx a,\psi\rx 1$ defines an $\mg$-linear map
$\ro:\,A\rightarrow X^\ast$, where $X^\ast$ is the $\mg$-dual of $X$.
We want to prove that $\ro$ is surjective. Since $X$ is 
finite-dimensional, it suffices to prove that $\ro(A)$ separates 
points on $X$. Let $0\neq\psi\in X$, and choose $c\in A$ with
$\lx c,\psi\rx\neq 0$. After decomposing  with respect to the bigrading of 
$D_\gh$, we may assume that
 $\lx c,\psi\rx\in (\D)_{\de\de}$ for some $\de$,
which implies $\ro(c)(\psi)=\lx c,\psi\rx T_\de\neq 0$.
This proves that $\ro$ is surjective. In particular, we can find
$b\in A$ with $\ro(b)(\psi_i)=\de_{i1}$. 
\end{proof}

\begin{lemma}\label{il}
Let $A$ be an $\gh$-coalgebroid.
Then the map 
$$\lx a\ot b,\iota(\phi\ot\psi)\rx=\lx a,\phi\rx\,T_\be\,
\lx b,\psi\rx,\qquad \phi\in A_{\al\be}',\ \psi\in A_{\be\ga}',
$$
defines an $\gh$-algebra embedding
$\iota:\, A'\wt\ot A'\hookrightarrow(A\wh\ot A)'.$
\end{lemma}

\begin{proof}
First we must check that $\iota$ is well-defined, that is, that
\begin{align*}\lx a\ot b,\iota(\mu_r(f)\phi\ot\psi)\rx
&=\lx a\ot b,\iota(\phi\ot\mu_l(f)\psi)\rx,\\
\lx a\mu_l(f)\ot b,\iota(\phi\ot\psi)\rx&=\lx a\ot\mu_l(f)b,
\iota(\phi\ot\psi)\rx,\\
\lx a\mu_r(f)\ot b,\iota(\phi\ot\psi)\rx&=\lx a\ot\mu_r(f)b,
\iota(\phi\ot\psi)\rx
\end{align*}
for $f\in\mg$. This is straight-forward. For instance, to prove the second
identity, we assume $\phi\in A'_{\al\be}$, $a\in A_{\ga\de}$. We must prove 
that 
$$\lx a\mu_l(f),\phi\rx\,T_\be\,\lx b,\psi\rx=\lx a,\phi\rx\,T_\be\,
\lx\mu_l(f)b, \psi\rx,$$
or equivalently, 
by \eqref{dr} and \eqref{ap0}, that
$$T_{-\ga} f\circ\lx a,\phi\rx\circ T_\be\circ\lx b,\psi\rx
=\lx a,\phi\rx\circ T_\be\circ f\circ \lx b,\psi\rx.$$
This is clear from the fact that 
$\lx a,\phi\rx\in (D_\gh)_{\ga+\be,\de+\al}=\mg T_{-\ga-\be}$.

It is  easy to check that
$\iota$ maps into the subspace $(A\wh\ot A)'$ of $\Hom_{\C}(A\wh\ot A,D_\gh)$
and that $\iota$ is an $\gh$-prealgebra homomorphism. 
Let us write out the
proof that $\iota$ is multiplicative, that is, that
\begin{equation}\label{im}\left\lx a\ot b,
\iota\big((\phi_1\ot\phi_2)(\psi_1\ot\psi_2)\big)\right\rx=\big\lx a\ot b,
\iota(\phi_1\ot\phi_2)\iota(\psi_1\ot\psi_2)\big\rx.
\end{equation}
Assume that $\phi_1\in A'_{\bullet\al}:=\bigoplus_\ep A'_{\ep\al}$, 
$\psi_1\in  A'_{\bullet\be}$, and  write
$\De(a)=\sum_i a_i'\ot a_i''$, $\De(b)=\sum_j b_j'\ot b_j''$,
where $a_i'\in A_{\bullet\ga_i}$,
$b_j'\in A_{\bullet\de_j}$. The left-hand side
of \eqref{im} is then
$$\lx a,\phi_1\psi_1\rx\,T_{\al+\be}\,\lx b,\phi_2\psi_2\rx=\sum_{ij} \lx a_i'
,\phi_1\rx\,T_{\ga_i}\,\lx a_i'',\psi_1\rx\,T_{\al+\be}\,\lx b_j',\phi_2\rx
\,T_{\de_j}\,\lx b_j'',\psi_2\rx,$$
while the right-hand side is
\begin{equation*}\begin{split}&
\sum_{ij}\lx a_i'\ot b_j',\iota(\phi_1\ot\phi_2)\rx\,T_{\ga_i+\de_j}\,
\lx a_i''\ot b_j'',\iota(\psi_1\ot\psi_2)\rx\\
&\quad=\sum_{ij}\lx a_i',\phi_1\rx\,T_\al\,\lx b_j',\phi_2\rx\,
T_{\ga_i+\de_j}\,\lx a_i'',\psi_1\rx\,T_{\be}\,\lx b_j'',\psi_2\rx.
 \end{split}\end{equation*}
These expressions are equal since
$\lx a_i'',\psi_1\rx\in\mg\, T_{-\ga_i-\be}$, 
$\lx b_j',\phi_2\rx\in\mg\, T_{-\de_j-\al}$.

To prove that $\iota$ is injective, let
 $0\neq x\in A'\wt\ot A'$. 
 We choose a representative $\sum_{i=1}^n\phi_i\ot\psi_i$ 
of $x$ such that 
$(\psi_i)_{i=1}^n$ are  independent over $\mu_l(\mg)$
and $\phi_i\neq 0$ for all $i$. Choosing $b$ as in Lemma \ref{ld},
we have 
 $\lx a\ot b,\iota(x)\rx 1=\lx a,\phi_1\rx 1$. Again by Lemma \ref{ld},
there exists $a\in A$ with $\lx a,\phi_1\rx 1\neq 0$. Thus $\iota(x)\neq 0$.
This completes
the proof.
\end{proof}

We now let $A$ be an $\gh$-bialgebroid, and 
 apply Lemma \ref{dhl} to the multiplication, or rather to its 
quotient $m:\,A\wh\ot A\rightarrow A$. This gives a map
$m':\,A'\rightarrow(A\wh\ot A)'$ defined by
$\lx x\ot y,m'(\phi)\rx=\lx xy,\phi\rx.$
Now let
$$A^\ast=\left\{ \phi\in A'\ ;\ m'(\phi)\in\Ima(\iota)\right\}.$$
It follows from the fact that $m'$ and $\iota$ are $\gh$-algebra homomorphisms
that $A^\ast$ is an $\gh$-subalgebra of $A'$.

\begin{lemma}The map  $\De=\iota^{-1}\circ m'\big|_{A^\ast}$ is an 
$\gh$-algebra homomorphism $A^\ast\rightarrow A^\ast\wt\ot A^\ast$.
\end{lemma}

\begin{proof}
It remains to prove that $\De$ takes values in  $A^\ast\wt\ot A^\ast$.
Since $\De$ preserves the bigrading, 
it is enough to consider $\De(\chi)$ for $\chi\in A^\ast_{\al\be}$.
We choose a representative  $\sum_{i=1}^n\phi_i\ot\psi_i$ of $\De(\chi)$
such that $\phi_i\in A'_{\al\ga_i}$, $\psi_i\in A'_{\ga_i\be}$,
$(\phi_i)_{i=1}^n$ are linearly independent over $\mu_r^{A'}(\mg)$
and  $(\psi_i)_{i=1}^n$  are independent over $\mu_l^{A'}(\mg)$.
Then
$$\lx ab,\chi\rx=\sum_{i=1}^n \lx a,\phi_i\rx\,T_{\ga_i}\,\lx b,\psi_i\rx,
\qquad a,\, b\in A. $$
Choosing $b$ as in Lemma \ref{ld} and applying the above identity
to $(ac)b=a(cb)$ gives
$$\lx ac,\phi_1\rx 1=\sum_{i=1}^n\lx a,\phi_i\rx\,T_{\ga_i}\,\lx cb,\psi_i\rx
1, \qquad a,\,c\in A.$$
We may assume that $b=\sum_\ep b_\ep$ with $b_\ep\in A_{\ep+\ga_1-\be,\ep}$.
Define $\tilde \psi_i:\,A\rightarrow \D$ through 
 $\lx c,\tilde \psi_i\rx=\sum_\ep\lx cb_\ep,\psi_i\rx T_{\ep}$. 
It is then easily checked that
$\tilde\psi_i\in A_{\ga_i\ga_1}'$, which gives
$$m'(\phi_1)=\iota\left(\sum_{i=1}^n\phi_i\ot\tilde\psi_i\right)$$
so that $\phi_1\in A^\ast$. By symmetry, $\phi_i\in A^\ast$ for all $i$,
and by a similar argument we may conclude that also $\psi_i\in A^\ast$.
This completes the proof.
\end{proof}

It is now clear how to define the dual of an $\gh$-bialgebroid.

\begin{definition} For $A$ an $\gh$-bialgebroid, we define the dual
$\gh$-bialgebroid $A^\ast$ to be the $\gh$-subalgebra $A^\ast\subseteq A'$,
equipped with the coproduct $\De=\iota^{-1}\circ m'$ and the counit
$\ep:\,A^\ast\rightarrow D_{\gha}$ defined by $\ep(\psi)=\lx 1,\psi\rx$.
\end{definition}

It is easy to check that the coproduct and counit axioms are satisfied.

\begin{proposition}
If $A$ is an $\gh$-Hopf algebroid, then $A^\ast$ is an $\gh$-Hopf algebroid
with the antipode $S^{A^\ast}(\phi)=S^{D_\gh}\circ\phi\circ S^A$.
\end{proposition}

\begin{proof}
It is easy to check that $S(A')\subseteq A'$. To see that $S$ preserves 
$A^\ast$, choose $\phi\in A^\ast$ with $\De(\phi)=\sum_i\phi_i'\ot\phi_i''$,
and check that $m'(S(\phi))=\iota\left(\sum_iS(\phi_i'')\ot S(\phi_i')\right)$.
Next,  check that $S(A_{\al\be}^\ast)\subseteq A^\ast_{-\be,-\al}$.
The condition \eqref{sf} is then equivalent to 
$S(\mu_l(f)\phi)=S(\phi)\mu_r(f)$,
$S(\phi\mu_r(f))=\mu_l(f)S(\phi)$, which is easily verified.

To check the first identity in \eqref{anti}, we let
$\phi\in A^\ast_{\al\be}$ and $x\in A_{\ga\de}$ with
$\De(\phi)=\sum_i\phi_i'\ot\phi_i''$, $\phi_i''\in A_{\xi_i\be}^\ast$,
$\De(x)=\sum_j x_j'\ot x_j''$, $x_j''\in A_{\eta_j\de}$. Then
\begin{equation*}\begin{split}&
\lx x,m(\id\ot S)\De(\phi)\rx=\sum_i\lx x,\phi_i'\,S(\phi_i'')\rx
=\sum_{ij}\lx x_j',\phi_i'\rx\,T_{\eta_j}\,\lx x_j'',S(\phi_i'')\rx\\
&\quad=\sum_{ij}\lx x_j',\phi_i'\rx\,T_{\eta_j}\,S^{D_\gh}\!
\left(\lx S(x_j''),\phi_i''\rx\right)
=\sum_{ij}\lx x_j',\phi_i'\rx\,T_{\xi_i}\,\lx S(x_j''),\phi_i''\rx\,
T_{\be-\de},\end{split}\end{equation*}
where we used that $\lx S(x_j''),\phi_i''\rx\in(D_{\gh})_{\be-\de,\xi_i-\eta_j}$.
By the definition of $\De^{A^\ast}$ and by \eqref{anti} for $S^A$, this equals
$$\sum_j \lx x_j'S(x_j''),\phi\rx\, T_{\be-\de}
=\lx \mu_l(\lx x,\ep\rx 1),\phi\rx\, T_{\be-\de}
=\lx x,\ep\rx 1\circ\lx 1,\phi\rx\circ T_{\be-\de}.$$
On the other hand,
$$\lx x,\mu_l(\ep(\phi)1)\rx=\lx x,\mu_l(\lx 1,\phi\rx 1)\rx
=\lx 1,\phi\rx 1\circ \lx x,\ep\rx.$$
These two expressions are equal since
 $\lx 1,\phi\rx\in\mg T_{-\be}$, $\lx x,\ep\rx\in\mg T_{-\de}$.
This proves the first identity in \eqref{anti} and
the second one may be proved similarly.
\end{proof}

\begin{proposition}\label{kap}
For $A$ an $\gh$-Hopf algebroid, there exists a homomorphism
of $\gh$-Hopf algebroids 
$\ka:\,A\rightarrow A^{\ast\ast}$, defined by 
$\lx \phi,\ka(x)\rx=\lx x,\phi\rx$.
\end{proposition}

The proof is straight-forward. 

\begin{example}\label{die}
Since  $\D$ is the unit object for $\gh$-algebras (as a tensor category)
and  $\I$ the unit object for $\gh$-coalgebroids,
it is natural to expect that $\I$ and $\D$ are mutually dual. 
Indeed, it is easy to check that 
 $\I'\simeq \D$ through the isomorphism
\begin{equation}\label{ip}
\lx f\ot g,x\rx=f\circ x\circ g, \qquad f,\,g\in\mg,\ x\in\D.\end{equation}
Moreover, the embedding
 $\iota:\,\D\wt\ot \D\hookrightarrow (\I\wh\ot \I)'$
of Lemma \ref{il} is an isomorphism. In fact, it can be
 identified with the identity map on $\D$ through the canonical
isomorphisms \eqref{can}, \eqref{can2}. Thus we have also an equivalence
of $\gh$-Hopf algebroids
$\I^\ast\simeq \D$.

Next, we let  $\phi\in\D'$  and introduce the function
$f_\phi:\,\gha\oplus\gha\rightarrow \C\cup\{\infty\}$ through
$$f_\phi(\la,\mu)=(\lx T_{\mu-\la},\phi\rx 1)(\la).$$
Then $\phi\mapsto f_\phi$ gives an equivalence 
$$\D'\simeq\left\{f:\,\gha\oplus\gha\rightarrow \C\cup\{\infty\};
\ \la\mapsto f(\la,\la+\mu)\ 
\text{is meromorphic for all}\ \mu\in\gha\right\}.$$
 The space $\D^\ast$ is identified with the subspace of
tensors:
$$\D^\ast
\simeq\left\{f\in \D';\ f(\la,\mu)=\textstyle\sum_{\text{finite}}
\displaystyle g_i(\la)h_i(\mu)\right\}.$$
We have set this up so that the map $\ka:\,\I\rightarrow \I^{\ast\ast}
=\D^{\ast}$ of Proposition \ref{kap}
is the natural embedding $\ka(g\ot h)(\la,\mu)=g(\la)h(\mu)$.
Strictly speaking  $\D^{\ast}$ is  bigger than $\I$,
as can be seen by considering  $f(\la,\mu)=g(\la)/g(\mu)$
with $g$  a discontinuous character, $g(\la+\mu)=g(\la)g(\mu)$.
\end{example}

\begin{example}
This is a pathological example which
shows that the map $\ka$ of Proposition \ref{kap} is not 
always an embedding. Let $f\in\mg$ and let $K_f\subseteq \I$ be 
the principal ideal generated by $f\ot 1-1\ot f$. One may check that  $\I/K_f$
is an 
$\gh$-Hopf algebroid (cf.~Lemma \ref{hq} below). As in the previous example,
any $\phi\in (\I/K_f)'$ is given by an element $x\in \D$ through
\eqref{ip}. However, $x\in \D$ defines an element of  $(\I/K_f)'$
if and only if $f\circ x=x\circ f$, which for generic $f$ implies 
$x\in(\D)_{00}$. But then $\phi$ vanishes on every element of the
form $a=g\ot 1-1\ot g$, $g\in\mg$, so that $\ka(a)=0$. Since not
necessarily $a\in K_f$, we have proved that $\ka$ is not injective. 
\end{example}

\subsection{Dual representations}
In this section we will establish a  correspondence
between corepresentations of an $\gh$-bialgebroid $A$ and $\gh$-representations
of the dual algebroid $A^\ast$.

When $V$ is an $\gh$-space we define $V^+$ to be the $\gh$-space
$V^+=\bigoplus_\al V_\al^+$, where $V_\al^+=\Hom_{\mg}(V_{-\al},\mg)$
(note the minus sign!),
 with the $\mg$-linear structure
$$\langle v,\xi f\rangle=\,\lx v,\xi\rx f,\qquad v\in V,\ 
\xi\in V^+,\ f\in\mg, $$ 
where we write $\lx v,\xi\rx=\xi(v)\in\mg$. Note that, by definition,
$\lx fv,\xi\rx=f\lx v,\xi\rx$.

\begin{proposition}\label{drp}
Let $\pi:\,V\rightarrow A\wt\ot V$ be a corepresentation of an 
$\gh$-bialgebroid $A$ on an $\gh$-space $V$. 
Then there exists an $\gh$-representation
$\pi^+$ of $A^{\ast}$ on $V^+$, defined by 
$$\lx v,\pi^+(\phi)\xi\rx=\lx \pi(v),\phi\ot \xi\rx,\qquad
v\in V,\ \phi\in A^\ast,\ \xi\in V^+,$$
where the pairing on the right-hand side is defined by
$$\lx a\ot v,\phi\ot\xi\rx=\left(\lx a,\phi\rx \,T_{\be}\,\lx v,\xi\rx\right) 
1,\qquad a\in A_{\al\be}.$$
\end{proposition}

We omit the straight-forward proof, which is largely parallel to the first
part of the proof of Lemma \ref{il}. 

Let $\{v_k\}_k$
be a basis of $V$, $t_{kj}\in A$ the corresponding matrix elements
as in \eqref{cme} and $v_k^+\in V^+$  the dual basis elements defined by
 $\lx f v_j,v_k^+\rx=f\de_{jk}$. Then the representation $\pi^+$ 
is given by
\begin{equation}\label{drx}
\pi^+(\phi)v_k^+=\sum_j v_j^+\lx t_{jk},\phi\rx 1.
\end{equation}

The coproduct $\De$ defines a corepresentation of any 
 $\gh$-coalgebroid on itself (the regular corepresentation). It is less
 easy to define  the regular $\gh$-representation 
of an $\gh$-algebra on itself. The naive definition $\pi(a)b=ab$
does not work since then \eqref{md} would imply a connection between
the left and right moment map. However, we can define it using duality.

\begin{definition} For $A$ an $\gh$-algebra, we define the 
regular representation of $A$ to be the $\gh$-representation
$\ro\circ\ka:\,A\rightarrow D_{A^{\ast+}}$, where 
$\ro:\, A^{\ast\ast}\rightarrow D_{A^{\ast+}}$ is the 
dual of the regular 
corepresentation of $A^\ast$ and $\ka:\,A\rightarrow A^{\ast\ast}$
is as in \emph{Proposition \ref{kap}}.
\end{definition}

\begin{example}
Consider the regular representation of the $\gh$-algebra $I_\gh$. 
According to Example \ref{die}, it is realized on the $\gh$-space $D_\gh^+$.
We identify elements of $I_\gh$ with functions on $\gha\oplus\gha$ through
$(f\ot g)(\la,\mu)=f(\la)g(\mu)$ and  
elements $\xi\in\D^+$ with functions through
$$\xi(\la,\mu)=\lx T_{\mu-\la},\xi\rx(\la).$$
Then $\D^+$ is identified with the space of finite sums
$$f(\la,\mu)=\sum_{\al\in\gha}g_\al(\la)\,\de_{\mu-\la,\al},
\qquad g_\al\in\mg, $$
equipped with the $\gh$-space structure
$g(\la)\de_{\mu-\la,\al}\in(\D^+)_{\al\al}$, 
$$(g f)(\la,\mu)=f(\la) g(\la,\mu),\qquad 
g\in \D^+,\ f\in \mg.$$
 Note that it follows that
$(fg)(\la,\mu)=f(\mu)g(\la,\mu)$; this gives the connection between the
left and right moment maps needed to define the regular representation. 
It is the $\gh$-representation $\pi$ of $\I$ 
on $\D^+$  given by pointwise multiplication of functions: 
$$(\pi(g) h)(\la,\mu)=g(\la,\mu)h(\la,\mu),\qquad g\in\I,\ h\in\D^+.$$
\end{example}

For completeness, we state the following fact, which shows that 
the the duality of Proposition \ref{drp} extends to the level of tensor 
categories. Again we omit the straight-forward proof.

\begin{proposition}\label{til}
Let $V$ and $W$ be corepresentations of an $\gh$-bialgebroid $A$.
If $\Phi:\,V\rightarrow W$ is an intertwiner of corepresentations,
then $\lx v,\Phi^+(\xi)\rx=\lx \Phi(v),\xi\rx$ defines an intertwiner
of $\gh$-representations $\Phi^+:\,W^+\rightarrow V^+$. 

Moreover,
there exists an intertwiner (an equivalence in the finite-dimensional case)
of $\gh$-representations $\iota:\,W^+\wh\ot V^+\rightarrow (V\wh\ot W)^+$ 
defined by
$$\lx x\ot y,\iota(\eta\ot \xi)\rx=\lx x,\xi\rx
\,T_{-\al}(\lx y,\eta\rx),\qquad x\in V_\al. $$
\end{proposition} 

\subsection{Pairings and cobraidings}
\label{sps}

In view of the  results of \S \ref{ssad}, the following definition of
a pairing of $\gh$-bialgebroids is natural.

\begin{definition}\label{pd}
For $A$ and $B$ $\gh$-bialgebroids, we define a pairing to be a 
$\C$-bilinear map $\langle\cdot,\cdot\rangle:\,A\times B\rightarrow D_{\gh}$,
with
\begin{gather}
\nonumber\langle A_{\al\be},B_{\ga\de}\rangle\subseteq
(D_\gh)_{\al+\de,\be+\ga},\\
\begin{split}\label{mmp}\langle\mu_l(f)a,b\rangle&=\langle a,\mu_l(f)b\rangle
=f\circ\langle a,b\rangle,\\
\langle a\mu_r(f),b\rangle&=\langle a,b\mu_r(f)\rangle
=\langle a,b\rangle\circ f,\end{split}\\
\begin{split}\label{ppp}
\langle ab,c\rangle&=\sum_i\langle a,c_i'\rangle 
\,T_{\be_i}\,\langle b,c_i''\rangle,
\qquad \De(c)=\sum_i c_i'\ot c_i'',\qquad c_i''\in  B_{\be_i\ga},\\
\langle a,bc\rangle&=\sum_i\langle a_i',b\rangle 
\,T_{\be_i}\,\langle a_i'',c\rangle,
\qquad \De(a)=\sum_i a_i'\ot a_i'',\qquad a_i''\in  A_{\be_i\ga},\end{split}\\
\nonumber\langle a,1\rangle=\ep(a),\qquad \langle 1,b\rangle =\ep(b).
\end{gather}
A pairing is said to be non-degenerate if $\langle a,B\rangle=0\Rightarrow a=0$
and $\langle A,b\rangle=0\Rightarrow b=0$.
\end{definition}

It is clear that a pairing defines homomorphisms
(embeddings in the non-degenerate case) of $\gh$-bialgebroids
$A\rightarrow B^\ast$ and $B\rightarrow A^\ast$. 

A pairing between $\gh$-Hopf algebroids should in addition  satisfy
\begin{equation}\label{hap}
\langle Sa,b\rangle=S^{D_\gh}(\langle a,Sb\rangle).\end{equation}

For $A$ an  $\gh$-bialgebroid, we denote by
$A^{\cop}$  the  $\gh$-bialgebroid which equals $A$ as
an associative algebra but has the opposite $\gh$-coalgebroid  structure, 
that is,
\begin{gather*} A_{\al\be}^{\cop}=A_{\be\al},\qquad \mu_l^{A^{\cop}}=\mu_r^A,
\qquad \mu_r^{A^{\cop}}=\mu_l^A,\qquad
\De^{A^{\cop}}=\si\circ\De^A,\qquad \ep^{A^{\cop}}=\ep^A,\end{gather*}
where $\si(x\ot y)=y\ot x$. If $A$ is an $\gh$-Hopf-algebroid
with invertible antipode $S$, then $A^{\cop}$ is an $\gh$-Hopf-algebroid
with antipode $S^{-1}$.

For Hopf algebras,
 pairings on $A^{\cop}\times A$ are of special interest
since they are related to quasitriangular or braided structures on $A^\ast$.
The following definition will turn out to be appropriate in the
dynamical case. 

\begin{definition}
A cobraiding on an $\gh$-bialgebroid $A$ is a 
pairing  $\langle\cdot,\cdot\rangle$ on  $A^{\cop}\times A$ which satisfies
\begin{equation}\label{cb}\sum_{ij}\mu_l(\lx a_i',b_j'\rx 1)\,a_i''b_j''
=\sum_{ij}\mu_r(\lx a_i'',b_j''\rx 1)\,b_j'a_i',\end{equation}
where $\De^A(a)=\sum_i a_i'\ot a_i''$ and similarly for $b$.
\end{definition}

Note that, in terms of the equivalence \eqref{can}, we may write \eqref{cb} as
$$ \sum_{ij} a_i''b_j''\ot \lx a_i',b_j'\rx =\sum_{ij}
\lx a_i'',b_j''\rx \ot b_j'a_i'. $$

To prove that something is a  cobraiding one needs the following lemma,
which shows that it is enough to verify \eqref{cb} for a set of generators.
We omit the straight-forward proof; cf. \cite{k} for the case of Hopf
algebras.

\begin{lemma}\label{cbvl}
If the condition \eqref{cb} holds with $(a,b)$ replaced by $(x,y)$, $(x,z)$
and $(y,z)$, then it also holds for $(xy,z)$ and $(x,yz)$. 
\end{lemma}

A cobraiding on an $\gh$-bialgebroid yields a braiding 
(in the dynamical sense) on its corepresentations. 

\begin{proposition}\label{cbp}
Let $A$ be an $\gh$-bialgebroid equipped with a cobraiding, and let $V$ and
$W$ be two corepresentations of $A$. Choosing bases and introducing matrix 
elements as in \eqref{cme}, the equation
$$\Phi(v_j\ot w_k)=\sum_{lm}\lx t_{kl}^W,t_{jm}^V\rx 1\,w_l\ot v_m$$
defines an intertwiner $\Phi:\,V \wh\ot W\rightarrow W \wh\ot V$.

 Let $V_0\subseteq V$, $W_0\subseteq W$
be the complex subspaces spanned by $\{v_k\}_k$, $\{w_k\}_k$, and define
$\Phi_0:\,\gha\rightarrow\Hom_\gh(V_0\ot W_0,W_0\ot V_0)$ through
 $$\Phi_0(\la)(v_j\ot w_k)=\sum_{lm}(\lx t_{kl}^W,t_{jm}^V\rx 1)(\la)
\,w_l\ot v_m.$$
Define $R=R_{V_0W_0}:\,\gha\rightarrow
\End_\gh(V_0\ot W_0)$ by
 $\Phi_0=\si\circ R$ where $\si(x\ot y)=y\ot x$. 
Then, given also 
a third corepresentation $U$ with distinguished complex subspace
$U_0$, $R$  satisfies the \emph{QDYB} equation of the form
\begin{equation}\label{uvw}\begin{split}&
R^{12}_{U_0V_0}(\la-h^{(3)})\,R^{13}_{U_0W_0}(\la)\,
R^{23}_{V_0W_0}(\la-h^{(1)})\\
&\quad=R^{23}_{V_0W_0}(\la)\,
R^{13}_{U_0W_0}(\la-h^{(2)})\,R^{12}_{U_0V_0}(\la).\end{split}\end{equation}
\end{proposition}

Note that it is not possible to factor $\Phi$ as $\si\circ \tilde R$ for some
$\tilde R\in\End_\gh(V\wh\ot W)$, since
$\si$ does not make sense as an operator on $V\wh\ot W$. 
As was pointed out in \cite{ev},
it is  this non-naturality of the flip map that
gives rise to the  QDYB equation instead of the elementary braid relation
$R^{12}R^{13}R^{23}=R^{23}R^{13}R^{12}$.

\begin{proof}
It follows from \eqref{ix} and \eqref{tpx} that the equation
$$\Phi(v_j\ot w_k)=\sum_{lm}\Phi_{jk}^{lm}\,w_l\ot v_m $$
defines an intertwiner if and only if
\begin{equation}\label{mei}\sum_{ij}\mu_l(\Phi_{kl}^{ij})t_{im}^W\,t_{jn}^V
=\sum_{ij}\mu_r(\Phi_{ji}^{mn})
t_{kj}^V\, t_{li}^W. \end{equation}
Choosing $a=t_{lm}^W$ and $b=t_{kn}^V$ in \eqref{cb} and
using \eqref{dt}, we see that this indeed holds for
$\Phi_{jk}^{lm}=\lx t_{kl}^W,t_{jm}^V\rx 1$.

To prove the second statement we apply 
$x\mapsto\lx x,t_{pq}^U\rx$ to both sides of \eqref{mei}.
Writing  $R^{jk}_{lm}=\Phi^{ml}_{jk}=\lx t_{km}^W,t_{jl}^V\rx 1$, the
 left-hand side gives
\begin{equation*}\begin{split}
&\sum_{ij}\big\lx t_{im}^W\,t_{jn}^V\,\mu_r^{A^{\cop}}(T_{\om(i)+\om(j)}
R^{kl}_{ji}),t_{pq}^U \big\rx\\
&\quad=\sum_{hij}\lx t_{im}^W,t_{ph}^U \rx\circ T_{\om(h)}
\circ \lx t_{jn}^V,t_{hq}^U\rx\circ T_{\om(i)+\om(j)}\circ R^{kl}_{ji}
\circ T_{-\om(i)-\om(j)}\\
&\quad=\sum_{hij}R^{pi}_{hm}\circ T_{-\om(m)}\circ R^{hj}_{qn}\circ 
T_{\om(i)-\om(h)}\circ R^{kl}_{ji}\circ T_{-\om(i)-\om(j)}\\
&\quad=\sum_{hij}R^{pi}_{hm}(\la)\,R^{hj}_{qn}(\la-\om(m))
\,R^{kl}_{ji}(\la-\om(p))\, T_{-\om(m)-\om(n)-\om(q)}.
\end{split}\end{equation*}
Computing the right-hand side similarly gives the identity
$$\sum_{hij}R^{pi}_{hm}(\la)\,R^{hj}_{qn}(\la-\om(m))
\,R^{kl}_{ji}(\la-\om(p))=\sum_{hij}R^{ji}_{nm}(\la)\,R^{pk}_{hj}(\la)
\,R^{hl}_{qi}(\la-\om(j)). $$
Replacing $(h,i,j,k,l,m,n,p,q)\mapsto(x,z,y,b,c,f,e,a,d)$
and comparing with \eqref{dyx}, we see that this is indeed the QDYB equation.
\end{proof}

\subsection{Cobraidings on dynamical quantum groups}

We now turn to the case of an
 $\gh$-bialgebroid $A_R$
constructed via the generalized FRST
construction from a matrix $R$. 
We will show that the  QDYB equation for $R$ gives
 a cobraiding on  $A_R$.

\begin{proposition}\label{lp}
In the setting of \emph{\S \ref{ssf}}, let  
$L:\,\gh^\ast\rightarrow\End_\gh(V\ot V)$ be a
meromorphic function.
We introduce its matrix elements 
 $L^{ab}_{xy}\in\mg$ as in \eqref{me}\emph{;} note that
 $L^{ab}_{xy}=0$ if $\om(x)+\om(y)\neq \om(a)+\om(b)$.
Then the following three statements are equivalent\emph{:}

\begin{enumerate}[\emph(i\emph)]
\item There exists a pairing $A_R^{\cop}\times A_R\rightarrow D_\gh$ defined by
\begin{equation}\label{lpl}\langle L_{ij},L_{kl}\rangle=L_{ik}^{jl}
\,T_{-\om(i)-\om(k)}.\end{equation}
\item $L$ satisfies the  relations
\begin{subequations}\label{rl}\begin{align}\label{rl3}
R^{12}(\la-h^{(3)})L^{13}(\la)L^{23}(\la-h^{(1)})&=L^{23}(\la)
L^{13}(\la-h^{(2)})R^{12}(\la),\\
\label{rl4}L^{12}(\la-h^{(3)})L^{13}(\la)R^{23}(\la-h^{(1)})&=R^{23}(\la)
L^{13}(\la-h^{(2)})L^{12}(\la).\end{align}\end{subequations}
\item There exist two $\gh$-representations $\pi$, $\ro$ on  $\mg\ot V$
of $A_R$, $A_R^{\cop}$, 
respectively, 
 defined by
\begin{align*}\pi(L_{ij})(g\ot v_k)&=\sum_l L_{il}^{jk}\,T_{-\om(j)}g\ot v_l,\\
\ro(L_{ij})(g\ot v_k)&=\sum_l L_{ki}^{lj}\,T_{-\om(i)}g\ot v_l.
 \end{align*}
\end{enumerate}
If these conditions are satisfied, the pairing is a cobraiding if and
only if
\begin{equation}\label{cbc}\sum_{xy} \lx L_{ax},L_{cy}\rx(\la)\,L_{xb}L_{yd}=
\sum_{xy} \lx L_{xb},L_{yd}\rx(\mu)\,L_{cy}L_{ax}\qquad \text{for all }
a,\,b,\,c,\,d.\end{equation}
\end{proposition}

Choosing $L=R$ in Proposition \ref{lp} and recalling \eqref{rll}
leads to the following algebraic interpretation of the
QDYB equation.

\begin{corollary}\label{rp}
The following two conditions are equivalent\emph{:}
\begin{enumerate}[\emph(i\emph)]
\item There exists a pairing $A_R^{\cop}\times A_R\rightarrow D_\gh$ defined by
$$\langle L_{ij},L_{kl}\rangle
=R_{ik}^{jl}\,T_{-\om(i)-\om(k)}.$$
\item $R$ satisfies the \emph{QDYB} equation \eqref{dy}.
\end{enumerate}
Moreover, this pairing is automatically a cobraiding.
\end{corollary}

\begin{proof}[Proof of \emph{Proposition \ref{lp}}.]
 Let $X$ and $Y$ denote
the left- and right-hand sides of \eqref{rll}, respectively.
Assuming the existence of a pairing, we have 
\begin{equation*}\begin{split}
\langle L_{ef},X \rangle&=\sum_{xyz} R^{xy}_{ac}(\la)\,
\langle L_{zf},L_{xb} \rangle\, T_{\om(z)}\,\langle L_{ez},L_{yd} \rangle\\
&=\sum_{xyz}  R^{xy}_{ac}(\la)\, L^{fb}_{zx}(\la)\, L_{ey}^{zd}(\la-\om(x))
\,T_{-\om(a)-\om(c)-\om(e)}.\end{split}\end{equation*}
Computing $\langle L_{ef},Y\rangle $ similarly 
one obtains the identity
$$\sum_{xyz}  R^{xy}_{ac}(\la)\, L^{fb}_{zx}(\la)\, L_{ey}^{zd}(\la-\om(x))=
\sum_{xyz} L^{fy}_{zc}(\la)\,L^{zx}_{ea}(\la-\om(c))\,R_{xy}^{bd}(\la-\om(f)).
$$
Replacing $(a,b,c,d,e,f,x,y,z)\mapsto(e,b,f,c,d,a,y,z,x)$ and comparing
with \eqref{dyx} wee see that this is equivalent to \eqref{rl4}.
Similarly, the identity $\langle X, L_{ef}\rangle=\langle Y, L_{ef}\rangle $
implies \eqref{rl3}. Thus, condition  (\emph{i}) implies condition (\emph{ii}).
The converse is proved similarly: using \eqref{mmp} and 
\eqref{ppp} we may extend 
\eqref{lpl} to  the algebra generated by $f(\la)$, $g(\mu)$, $L_{xy}$
subject to relations \eqref{zr}. By the above argument, \eqref{rl}
guarantees that the resulting form factors through the relations \eqref{rll}.

The equivalence of (\emph{ii}) and (\emph{iii}) is contained in 
\cite{ev}, Proposition 4.5. 

The ``only if''-part of the final statement follows by choosing 
$a=L_{ab}$, $b=L_{cd}$ in \eqref{cb}. The ``if''-part then follows using
Lemma \ref{cbvl}. 
\end{proof}

\section{Self-duality of the $\mathrm{SL}(2)$ dynamical quantum group}

\subsection{A pairing on $\fr(\mathrm{SL}(2))$}
We now turn to the case of Proposition \ref{lp} when $R$ is given by
\eqref{r}. As in the non-dynamical case, the solution $L=R$  
has  the disadvantage that the corresponding pairing does not factor through
the determinant relation \eqref{ddr}. Instead, one must work with 
$L=q^{-1/2}R$, which clearly satisfies \eqref{rl} and \eqref{cbc}.
We will denote the corresponding pairing by $\{\cdot,\cdot\}$.
To see that it factors through \eqref{ddr}, it
is enough to check that $\{ c,x\}=\{ x,c\}=\ep(x)$ for 
$c=\al\de-qF(\la)\ga\be$ and $x=\al$, $\be$, $\ga$, $\de$. This
yields proves the following proposition.

\begin{proposition}\label{scb} There exists a cobraiding $\{\cdot,\cdot\}$ 
 on $\fr$ determined by
$$\{L_{ij},L_{kl}\}=q^{-\frac 12}R_{ik}^{jl}\,T_{-\om(i)-\om(k)}.$$ 
\end{proposition}

Explicitly, the cobraiding is defined by
\begin{equation*}\begin{split}
\{ \al,\al\}&=q^{\frac 12}\,T_{-2},\qquad
\{ \al,\de\} =q^{-\frac 12}, \qquad \{ \de,\al\}=
q^{-\frac 12}G(\la),\qquad
\{ \de,\de\}=q^{\frac 12}\,T_2 ,\\
\{ \be,\ga\}&=q^{-\frac 12}\frac{q^{-1}-q}{q^{2(\la+1)}-1},\qquad
\{ \ga,\be\}=q^{-\frac 12}\frac{q^{-1}-q}{q^{-2(\la+1)}-1};
\end{split}\end{equation*}
all other combinations of generators give zero.
Although we will not need it, we note that the antipode axiom \eqref{hap}
is satisfied.

For our purposes (especially to simplify the proof of Lemma \ref{tll} below),
it will be more convenient to work with a related pairing. 
It is straight-forward to check that 
\begin{equation}\label{psi}
\begin{split}\Psi(\al)&=q^{\frac 12(\la-\mu)}\al,\qquad 
\Psi(\be)=q^{-\frac 12(\la+\mu)}F(\la)\ga,\\
\Psi(\ga)&=q^{\frac 12(\la+\mu)}\frac{1}{F(\mu)}\,\be,
\qquad \Psi(\de)=q^{\frac 12(\mu-\la)}\frac {F(\la)}{F(\mu)}\,\de
\end{split}\end{equation}
defines an isomorphism of  $\gh$-Hopf algebroids 
$\Psi:\,\fr\rightarrow\fr^{\cop}$.
Therefore, we may construct a pairing  on  $\fr\times\fr$ 
as $\lx x,y\rx=\{\Psi(x),y\}$. Explicitly, 
\begin{equation}\begin{split}\label{sp}
\langle \al,\al\rangle&=T_{-2},\qquad
\langle \al,\de\rangle =\langle \de,\al\rangle=1,\qquad
\langle \de,\de\rangle=G(\la+1)\,T_{2} ,\\
\langle \be,\be\rangle&=\frac{q-q^{-1}}{q^{\la}-q^{-\la}},\qquad
\langle \ga,\ga\rangle=\frac{q^{-1}-q}{q^{\la+2}-q^{-\la-2}}.
\end{split}\end{equation}

It is not hard to obtain the above pairing  from a 
general Ansatz, without using Proposition \ref{lp}. Note also that it
 is symmetric. In fact, one can prove that
 this is the case for all pairings on $\fr\times\fr$.

In the formal limit $f(\la) T_{-\al}\rightarrow f(-\infty)$
we obtain the highly degenerate pairing
\begin{equation*}\begin{split}\langle \al,\al\rangle&
=\langle \de,\de\rangle=\langle \al,\de\rangle=\langle \de,\al\rangle=1,\qquad
\langle \be,\be\rangle=\lx\ga,\ga\rx=0\end{split}\end{equation*}
on  $\mathcal F_q(\mathrm{SL}(2))$.
 We will see in Theorem \ref{deg}
 that, although the pairing \eqref{sp} is also degenerate, its kernel
is so small that  $\fr$ is ``almost'' self-dual. 
 In particular, 
$\mathcal F_R(\mathrm{SL}(2))$ can be viewed as a deformation
both of the function algebra $\mathcal F(\mathrm{SL}(2))$ and of the
universal enveloping algebra $\mathcal U(\mathrm{sl}(2))$.

\subsection{$\fr(\mathrm{SL}(2))$ as a deformed enveloping algebra}

In this section we work out the formal limit from $\fr(\mathrm{SL}(2))$ to 
$\uq$ explicitly. To this end we introduce the
elements
\begin{equation}\label{uf}
X_+=q^{-1}\frac{q^{\la+1}-q^{-\la-1}}{q-q^{-1}}\,\be,\qquad
X_-=-q\frac{q^{\mu}-q^{-\mu}}{q-q^{-1}}\,\ga,\qquad 
K=q^{\frac 12(\la-\mu)}\end{equation}
of $\fr$.
Note that these elements degenerate in the non-dynamical limit
$\la,\,\mu\rightarrow-\infty$, but not in the rational limit $q\rightarrow 1$.

\begin{proposition}\label{ufp}
The elements $X_+$, $X_-$, $K$, $K^{-1}$ generate a subalgebra of
$\fr$ isomorphic to $\uq$, that is,
\begin{equation}\label{uqr}\begin{split}
&X_+X_--X_-X_+=\frac{K^2-K^{-2}}{q-q^{-1}},\\
&KX_{\pm}=q^{\pm 1}X_{\pm}K,\qquad KK^{-1}=K^{-1}K=1.\end{split}\end{equation}
Moreover, the pairing \eqref{sp} satisfies
\begin{equation*}\begin{split}
\left\langle X_+,\left(\begin{matrix}\al&\be\\\ga&\de\end{matrix}
\right)\right\rangle&=\left(\begin{matrix}0& {F(\la-1)}\\0&0\end{matrix}
\right),\qquad 
\left\langle X_-,\left(\begin{matrix}\al&\be\\\ga&\de\end{matrix}
\right)\right\rangle=\left(\begin{matrix}0&0\\1/F(\la)&0\end{matrix}
\right),\\
 \left\langle K^{\pm},\left(\begin{matrix}\al&\be\\\ga&\de\end{matrix}
\right)\right\rangle&=
\left(\begin{matrix}q^{\pm\frac12}T_{-1}&0\\0&q^{\mp\frac 12}T_1
\end{matrix}\right),\end{split}\end{equation*}
and the Casimir element 
$$C=\frac{q^{-1}K^2+qK^{-2}-2}{(q^{-1}-q)^2}+X_+X_- $$
is related to the element $\Xi$ defined in \eqref{xi} by
$$C=\frac{1}{(q^{-1}-q)^2}\,\big(\Xi-2\big) .$$
\end{proposition}

The first relation in \eqref{uqr} corresponds to
$$\be\ga-\frac{F(\la)}{F(\mu-1)}\,\ga\be=I(\la,\mu),$$ 
which follows from the defining relations of $\fr$, and the remaining 
statements are straight-forward to check.

Note that $\uq$ is not embedded into $\fr$ as a Hopf subalgebra.
For instance, the coproduct rule $\De(\be)=\al\ot\be+\be\ot\de$ 
can be written as
\begin{equation}\label{dxp}
\De(X_+)=\left(\frac{q^{\la+1}-q^{-\la-1}}{q^{\mu+1}-q^{-\mu-1}}
\,\al\right)\ot X_+ +
X_+\ot \de.\end{equation}
To see how  the Hopf structure on $\uq$ arises in the non-dynamical
limit, rewrite the relations of $\fr$ using the generators 
$X_{\pm}$, $K^\pm$ instead of $\be$, $\ga$, and then formally  let
$\la,\,\mu\rightarrow-\infty$. For instance \eqref{dxp} becomes
$$\De(X_+)\sim \left(\frac{q^{\la+1}}{q^{\mu+1}}\,\al\right)\ot X_+ +
X_+\ot \de=K^2\al\ot X_+ +X_+\ot \de.$$
Thus we are effectively rescaling $\be$ 
and $\ga$ in the limit, and also considering $\la-\mu$ as fixed.

This results in a Hopf algebra $H$ with generators
$X_+$, $X_-$, $K$, $K^{-1}$, $\al$, $\de$ satisfying relations
\eqref{uqr}. The remaining relations say that $\al$ and $\de$ are central
elements with $\al\de=\de\al=1$. The Hopf structure on $H$ is given by
\begin{gather*}\begin{split}
\De(\al)&=\al\ot\al,\qquad \De(\de)=\de\ot\de,\qquad 
\De(K^{\pm})=K^\pm\ot K^\pm,\\
\De(X_+)&=X_+\ot\de+K^2\al\ot X_+,\qquad \De(X_-)
=X_-\ot K^{-2}\al+\de\ot X_-,\end{split}\\
\begin{split}\ep(\al)&=\ep(\de)=\ep(K^\pm)=1,\qquad \ep(X_\pm)=0,\\
S(\al)&=\de,\qquad S(\de)=\al,\qquad S(K^\pm)=K^\mp,\\
S(X_+)&=-K^{-2}X_+,\qquad S(X_-)=-X_-K^2.\end{split}\end{gather*}
It is clear from these relations 
that $\al-1$ and $\de-1$ generate a Hopf ideal $I$
of $H$ and that $H/I\simeq\uq$ with the standard Hopf structure.

It is curious that  $\uq$ is  contained 
in $\fr$ as an associative algebra. It would be interesting to know
whether the corresponding statement is true for the dynamical quantum
groups constructed in \cite{ev,ev2}.

\subsection{Finite-dimensional representations of $\fr(\mathrm{SL}(2))$}

Let $\pi$ be a corepresentation of $\fr$ on an $\gh$-space $V$.
In Proposition \ref{drp} we defined the dual $\gh$-representation
$\pi^+$ of $\fr^\ast$ on the dual space $V^+$. Composing it with the
morphism $\fr\rightarrow \fr^\ast$ induced from the pairing 
 \eqref{sp}, we obtain an $\gh$-representation 
of $\fr$ on $V^+$ which we also denote $\pi^+$.

We are interested in the dual of the finite-dimensional corepresentations
studied in \cite{kr}. Namely,
let $V_N$ be the subspace of $\fr$ spanned by
$\{\mu_l(\mc)\ga^{N-k}\al^k\}_{k=0}^N$, viewed as an $\gh$-space through
$\ga^{N-k}\al^k\in(V_N)_{2k-N}$, $fv=\mu_l(f)v$. Then 
$\pi_N=\De\big|_{V_N}$ defines
a corepresentation of $\fr$ on $V_N$. Its matrix elements 
$t_{kj}^N\in\fr$ are defined by
$$\De(\ga^{N-k}\al^k)=\sum_{j=0}^N t_{kj}^N\otimes 
\ga^{N-j}\al^j.$$
In \cite{kr}, the expression 
\begin{equation}\label{mex}
\begin{split}t_{kj}^N&=\sum_{l=\max(0,\,j+k-N)}^{\min(j,\,k)}
\qb{N-k}{j-l}\qb{k}{l}q^{j(j+2k-N)+l(3l-3k-3j+N)}\\
&\quad\times\frac{(q^{2(j-N-\mu-1)};q^2)_{j-l}}
{(q^{2(j+k-l-N-\mu-1)};q^2)_{j-l}}\,\ga^{j-l}\de^{N-k-j+l}\al^l\be^{k-l}
\end{split}\end{equation}
was derived; here we use the standard $q$-notation \cite{gr}
\begin{equation*}(a;q^2)_k=
\prod_{j=0}^{k-1}(1-aq^{2j}),\qquad
\genfrac{[}{]}{0pt}{}{n}{k}_{q^2}=
\frac{(q^2;q^2)_n}{(q^2;q^2)_k(q^2;q^2)_{n-k}}.
\end{equation*}

Let $\{v_k^+\}_{k=0}^N$ be the basis of $V_N^+$ dual to $\ga^{N-k}\al^k$,
that is, $\lx \ga^{N-j}\al^j,v_k^+\rx=\de_{jk}$. 
By \eqref{drx}, the representation 
$\pi_N^+:\,A\rightarrow D_{V_N^+}$
is given in this basis by
\begin{equation}\label{vnd}
\pi_N^+(a)v_j^+=\sum_{k=0}^N v_k^+ \lx t_{kj}^N,a\rx 1.\end{equation} 
The following lemma gives the action of the generators.

\begin{lemma}\label{tll}
The pairing \eqref{sp} satisfies
\begin{align*}
\lx t_{kj}^N,\al\rx&=\de_{kj}\,T_{N-2k-1},\\
\lx t_{kj}^N,\be\rx&=\de_{k,j+1}\,\frac{[k]}{[\la+N-2k+1]}\,T_{N-2k+1},\\
\lx t_{kj}^N,\ga\rx&=\de_{k+1,j}\,\frac{[k-N]}{[\la+2]}\,T_{N-2k-1},\\
\lx t_{kj}^N,\de\rx&=\de_{kj}\,\frac{[\la-k+1]\,[\la+N-k+2]}
{[\la+2]\,[\la+N-2k+1]}\,T_{N-2k+1},\\
\lx t_{kj}^N,\Xi\rx&=\de_{kj}\,(q^{N+1}+q^{-N-1})\,T_{N-2k},
\end{align*}
where we temporarily write $[x]=(q^{x}-q^{-x})/(q-q^{-1})$.
\end{lemma}

\begin{proof}
Iterating \eqref{ppp} and recalling the notation \eqref{lg} we have in general
$$\lx x_1x_2\dotsm x_N,L_{ij}\rx=\sum_{k_1,\dots,k_{N-1}\in\{+,-\}}
\lx x_1,L_{ik_1}\rx\,T_{\om(k_1)}\,\lx x_2,L_{k_1k_2}\rx\,T_{\om(k_2)}\dotsm
\lx x_N,L_{k_nj}\rx$$
(recall that $\om({\pm})=\pm 1$). Using this rule and \eqref{mex} to compute
$\lx t_{kj}^N,\de\rx$ we get two non-zero terms, corresponding to
$j=k=l+1$, $k_1=\dots=k_{N-1}=+$ and $j=k=l$, $k_1=\dots=k_{N-1}=-$.
Thus, 
\begin{equation*}\begin{split}\lx t_{kj}^N,\de\rx&=
\de_{kj}\left(\Big\lx
[N-k][k]\,\ga\de^{N-k-1}\al^{k-1}\be\,\frac{[\mu+k+1]}{[\mu+k]}\,,
\de\Big\rx+\lx \de^{N-k}\al^k,\de\rx\right)\\
&=\de_{kj}\bigg([N-k][k]\lx \ga,\ga\rx\,T_1\,\big(\lx \de,\al\rx
\,T_1\big)^{N-k-1}\big(\lx\al,\al\rx\, T_1\big)^{k-1}
\lx\be,\be\rx\,\frac{[\la+k+1]}{[\la+k]}\\
&\qquad +\big(\lx\de,\de\rx\,T_{-1}\big)^{N-k}\big(\lx\al,\de\rx\,
T_{-1}\big)^k\,
T_1\bigg)\\
&=\de_{kj}\bigg([N-k][k]\,\frac{-1}{[\la+2]}\,T_{N-2k+1}\,
\frac{1}{[\la]}\,\frac{[\la+k+1]}{[\la+k]}
 +\big(G(\la+1)\,T_{1}\big)^{N-k}\,T_{1-k}\bigg).
\end{split}\end{equation*}
Using that $G(\la)=F(\la)/F(\la-1)$ and $F(\la)=q^{-1}[\la+2]/[\la+1]$, we have
$$\big(G(\la+1)\,T_{1}\big)^{N-k}=\frac{F(\la+N-k)}{F(\la)}\,T_{N-k}
=\frac{[\la+N-k+2]\,[\la+1]}{[\la+N-k+1]\,[\la+2]}\,T_{N-k}. $$
Plugging this into the previous expression gives
\begin{equation*}\begin{split}\lx t_{kj}^N,\de\rx
&=\de_{kj}\,\frac{[\la+N-k+2]}{[\la+N-k+1]\,[\la+2]\,[\la+N-2k+1]}\\
&\quad\times\big\{[\la+1]\,[\la+N-2k+1]-[k]\,[N-k]\big\}
\,T_{N-2k+1}.\end{split}\end{equation*}
The expression in brackets equals 
$[\la-k+1][\la+N-k+1]$, which proves the statement for $\de$.
For $\al$, $\be$, $\ga$ the proof is similar but simpler, since we only
get one non-zero term. Finally we prove the statement for $\Xi$ using 
\eqref{dt}.  
\end{proof}

\begin{remark}\label{rc}
The $\gh$-representations obtained above may be compared  
with those studied in \cite{kr}.
Let, for $\om\in\mathbb C$, $\mathcal H^\om$ be the $\gh$-space with basis
$\{e_k\}_{k=0}^\infty$, where $e_k\in\mathcal H^\om_{\om+2k}$. In \cite{kr}
it was shown that the following equations define an $\gh$-representation
$\pi^\om$ of $\fr$ on $\mathcal H^\om$:
\begin{equation*}\begin{split}
\pi^\om(\al)e_k&= q^{-k}\frac{1-q^{2(\la-\om-k+1)}}{1-q^{2(\la-\om-2k+1)}}\,
e_k,\\
\pi^\om(\be)e_k&= \frac{(1-q^{2k})(1-q^{2(\om+k-1)})}{(1-q^{2(\la+1)})
(1-q^{2(\om+2k-\la-3)})}\,e_{k-1},\\
\pi^\om(\ga)e_k&=-q^{-1}e_{k+1},\\
\pi^\om(\de)e_k&=q^k \frac{1-q^{2(\la+1-k)}}{1-q^{2(\la+1)}}\,e_k.
\end{split}\end{equation*}
We will consider  the case $\om=-N\in\mathbb Z_{\leq 0}$, for which
$\mathcal H^\om_N=\bigoplus_{k=N+1}^\infty\mc e_k$ is
 an  $\gh$-submodule. We  write
$W_N=\mathcal H^\om/\mathcal H^\om_N$; this is a finite-dimensional 
$\gh$-representation of $\fr$. It is natural to expect that it is equivalent
to the representation $\pi_N^+$ defined in \eqref{vnd}. Indeed, one  
may check that this is true, the equivalence being given by
 $e_k=v_{N-k}^+N_k$, where 
$N_k(\la)= q^{k(\la+2)}(q^2;q^2)_k/(q^{2(\la+2)};q^2)_k.$
This  explains the fact (cf.~\cite{kr})
that the Clebsch--Gordan coefficients for  $\pi_N$ can be formally 
obtained from those of $\pi^\om$ by substituting $\om=-N$.
\end{remark}

\subsection{The radical of the pairing}\label{rps}

In this section we will compute the radical of 
  the pairing \eqref{sp}. We need the basic facts about 
ideals of $\gh$-Hopf algebroids;  these have not previously appeared
in the literature.

\begin{definition}\label{id}
A subspace $I$ of an $\gh$-Hopf algebroid $A$ is called an $\gh$-Hopf
ideal if the following conditions are satisfied\emph{:}
\begin{enumerate}[\emph(i\emph)]
\item $I$ is a two-sided ideal of $A$ as an associative algebra.
\item $(A_{\al\be}+I)\cap(A_{\ga\de}+I)=I$ if $(\al,\be)\neq(\ga,\de)$.
\item $\De(I)\subseteq A\wt\ot I+I\wt\ot A$.
\item $\ep(I)=0$ and $S(I)=I$.
\end{enumerate}
\end{definition}

Note that (\emph{ii}) guarantees that $A/I$ is bigraded.

\begin{lemma}\label{hq}
For $\Phi:\,A\rightarrow B$ a homomorphism of $\gh$-Hopf algebroids,
$\Ker\Phi=\{x\in A;\,\Phi(x)=0\}$ is an $\gh$-Hopf ideal. Moreover,
for any $\gh$-Hopf ideal $I$ in $A$, $A/I$ has a structure of
an $\gh$-Hopf algebroid such that the projection $A\rightarrow A/I$ is
an $\gh$-Hopf algebroid homomorphism.
\end{lemma}

The least trivial point is to prove that $I=\Ker\Phi$ satisfies
condition  (\emph{iii}) of Definition \ref{id}. Let us write
$\mu_c$ for the ``central'' moment map 
$\mu_c(f)(a\ot b)=\mu_r(f)a\ot b=a\ot\mu_l(f)b$
and consider each copy of $A$ in $A\wt\ot A$ as a vector space over
$\mu_c(\mg)$. By elementary linear algebra, we can then for each $x\in A$ 
write $\De(x)=\sum_{i,j\in\La_1\cup\La_2}\mu_c(f_{ij})x_i'\ot x_j''$,
where $x_i'$, $x_i''\in\Ker(\Phi)$ for $i\in\La_1$ and where the
families $(\Phi(x_i'))_{i\in\La_2}$, $(\Phi(x_i''))_{i\in\La_2}$  
are linearly independent over $\mu_c(\mg)$.
If $x\in\Ker(\Phi)$ we have
$$0=\De(\Phi(x))=(\Phi\ot\Phi)\De(x)
=\sum_{i,j\in\La_2}\mu_c(f_{ij})\,\Phi(x_i')\ot \Phi(x_j'').$$
Since $\wt\ot$  means tensor product over
$\mu_c(\mg)$, the linear independence
 gives $f_{ij}=0$ for  $i,\,j\in\La_2$, which implies
$\De(x)\in \Ker(\Phi)\wt\ot A+A\wt\ot\Ker(\Phi)$.

We want to compute the radical $\fr^\perp$ of the pairing
 \eqref{sp}, which is the
 space of  $x\in\fr$  such that $\lx x,y\rx=0$ for all $y\in\fr$,
or equivalently the kernel of the $\gh$-Hopf algebroid homomorphism
 $\Phi:\,\fr\rightarrow \fr^\ast$ induced by the pairing.

\begin{theorem}\label{deg}
Let $I\subseteq \fr$ be the left ideal generated by all 
elements $f(\la,\mu)$, where $f\in\mc\ot\mc$ satisfies $f(\la,\la+k)=0$
for all $k\in\Z$. 
Then  $\fr^\perp= I$. In particular, $I$ is an $\gh$-Hopf ideal and
$\fr/I$ is an $\gh$-Hopf algebroid with a 
 non-degenerate pairing  on $\fr/I\times\fr/I$. 
\end{theorem}

To get a better understanding of  $\fr/I$ we consider
its representations.

\begin{proposition}
Let $\pi$ be an $\gh$-representation of $\fr$ on an $\gh$-space $V$.
Then the following are equivalent\emph{:}
\begin{enumerate}[\emph(i\emph)]
\item $\pi$ factors to an $\gh$-representation of $\fr/I$.
\item  $V$ has the grading $V=\bigoplus_{k\in\mathbb Z} V_k$
\emph{(}where some $V_k$ may be zero\emph{)}.
\item Writing $K=q^{\frac 12(\la-\mu)}$ as in  \eqref{uf}, the spectrum of
$K$ is contained in $q^{\frac 12\mathbb Z}$.
\end{enumerate}
\end{proposition}

This follows immediately from the definitions. Note that  
condition (\emph {iii}) is natural both in
mathematics and physics:  it corresponds to particles
with half-integer spin and to representations of $\mathrm{SL}(2)$
rather than a covering group. 
We conclude that $\fr/I$  is  itself a good analogue of the
$\mathrm{SL}(2)$ Lie group. In particular, we can interpret Theorem
\ref{deg} as saying that ``the'' dynamical $\mathrm{SL}(2)$ quantum
group is self-dual.

Our main tool to prove  Theorem \ref{deg} is  the Peter--Weyl
theorem for $\fr$, proved in \cite{kr}, which says that the matrix
elements $t_{kj}^N$ form a basis for $\fr$ over $\mu_l(\mc)\mu_r(\mc)$.
 Roughly speaking, we will prove
Theorem \ref{deg} by constructing the dual basis with respect to
our pairing. For this we also need Lemma \ref{tll}.

First we construct the (truncated) projectors onto the isotypic 
components of the Peter--Weyl decomposition.

\begin{lemma}\label{pwp}
Let, for $L$, $M\in\Zp$, $P_{ML}\in\fr$ be the central element
$$P_{ML}=\prod_{\stackrel{0\leq l\leq L}{l\neq M}}(1-q^{l+1}\Xi+q^{2l+2})$$
and let $x$, $y\in\fr$. Assume that $x=\sum_{N\leq L} x_N$ with
$x_N=\sum_{kj}f_{kj}(\la,\mu)t_{kj}^N$. Then
$$\lx x,y P_{ML}\rx=C\lx x_M,y\rx, \qquad 0\neq C\in\C.$$
In particular, if $x\in\fr^\perp$, then each $x_N\in\fr^\perp$. 
\end{lemma}

\begin{proof}
Iterating \eqref{dt} gives in general
\begin{equation}\label{tx}\lx t_{kj}^N,x_1x_2\dotsm x_m\rx
=\sum_{l_1,\dots,l_{m-1}=0}^N
\lx t_{kl_1}^N,x_1\rx\,T_{2l_1-N}\,\lx t_{l_1l_2}^N,x_2\rx\,T_{2l_2-N}\dotsm
\lx t_{l_{m-1}j}^N,x_m\rx.\end{equation}
Using this and Lemma \ref{tll} we see that
$$\lx t_{kj}^N,p(\Xi)\rx=\de_{kj}\,p(q^{N+1}+q^{-N-1})\,T_{N-2k}$$\
for any complex polynomial $p$. 
Thus, for $N\leq L$,
\begin{equation*}\begin{split}\lx t_{kj}^N,P_{ML}\rx&
=\de_{kj}\prod_{\stackrel{0\leq l\leq L}{l\neq M}}(1-q^{l+1}(q^{N+1}+q^{-N-1})
+q^{2l+2})\,T_{N-2k}\\
&=\de_{kj}\prod_{\stackrel{0\leq l\leq L}{l\neq M}}(q^{N+1}-q^{l+1})
(q^{-N-1}-q^{l+1})\,T_{N-2k}
=\de_{kj}\,\de_{MN}\,C\,T_{N-2k}\end{split}\end{equation*}
with $C\neq 0$. This gives
$$\lx t_{kj}^N,yP_{ML}\rx=\sum_l\lx t_{kl}^N,y\rx\,T_{2l-N}\,
\lx t_{lj}^N,P_{ML}\rx=\de_{MN}\,C\,\lx t_{kj}^N,y\rx,$$
from which the statement readily follows. 
\end{proof}

Next we consider the projection from an isotypic component onto
the span of a single matrix element.

\begin{lemma}\label{mp}
For $y\in(\fr)_{st}$ and $f\in\mc\ot\mc$,
$$\big\lx f(\la,\mu)t_{kj}^N,\ga^{N-l}\be^Ny\ga^N\be^{N-m}\big\rx
=\de_{kl}\,\de_{jm}\,f(\la,\la+j-k-s)\,C(\la)\,
T_{-k}\,\lx t_{00}^N,y\rx\,T_{-j}$$
with $0\neq C\in\mc$.
\end{lemma}

To prove this we again use the expression \eqref{tx} together with Lemma 
\ref{tll}. To get a non-zero contribution
we must choose $l_1=k+1$, $l_2=k+2$ and so on up to $l_{N-l}=k+N-l$, which
gives $k\leq l$. The next $N$ indices decrease down to $l_{2N-l}=k-l$,
which gives $k\geq l$; this accounts for the factor $\de_{kl}$. Similarly,
starting from the right accounts for the factor $\de_{jm}$. 
Keeping track of the bigrading completes the proof.

We are now ready to prove Theorem \ref{deg}.
Assume that $x\in\fr^\perp$. Using first Lemma \ref{pwp} and then
 Lemma \ref{mp} we may assume that $x=f(\la,\mu)t_{kj}^N$, where
$$f(\la,\la+j-k-s)\lx t_{00}^N,y\rx= 0,\qquad y\in(\fr)_{st}.$$
Choosing $y=\al^s$ for $s\geq 0$ and $y=\de^{-s}$ for $s<0$, 
it follows from \eqref{tx} and Lemma \ref{tll} that
$\lx t_{00}^N,y\rx\neq 0$. Thus, $f(\la,\la+k)=0$ for all $k\in\Z$.
This shows that $\fr^\perp\subseteq I$. The reverse inclusion is easy to
prove.

\begin{remark}
The parameter $L$ in Lemma \ref{pwp} is only needed to ensure that
the projector belongs to $\fr$. In particular, we may consider
the  infinite product $P_{0\infty}$
as the projector from the regular to the trivial corepresentation, that is, 
as the  Haar
functional. More precisely, let $I_\C/I$ be the space $\mc\ot\mc$ modulo
the ideal generated by elements satisfying $f(\la,\la+k)=0$ for all $k\in\Z$.
We write $D_{\Z}=(I_\C/I)^\ast$; it is easy to check that
 $D_{\Z}\simeq \bigoplus_{k\in\Z}(D_\C)_{kk}$.
We define the Haar functional  on $\fr/I$ to be the map 
$h:\,\fr/I\rightarrow I_\C/I$ given by
$$h(f(\la,\mu)\,t_{kj}^N)=\de_{N0}\,f(\la,\mu), $$
cf.~also \cite{kr}, and its dual $h^\ast:\,D_{\Z}\rightarrow 
(\fr/I)^\ast$  by $\lx h(a),x\rx=\lx a,h^\ast(x)\rx$. 
It follows from Lemma \ref{pwp} that
$$h^\ast(fT_k)=\begin{cases} f(\la)\,\de^k P & k\geq 0,\\
 f(\la)\,\al^{-k}P, & k< 0,
\end{cases}$$
where
$$P=\prod_{l=1}^\infty
\frac{1-q^{l+1}\Xi+q^{2l+2}}{1-q^{l+1}(q+q^{-1})+q^{2l+2}}
=\frac{(q^2\xi,q^2\xi^{-1};q)_\infty}{(q,q^3;q)_\infty}; $$
here $\xi$ is a formal quantity satisfying $\xi+\xi^{-1}=\Xi$
and we use the standard notation
$$(a,b;q)_\infty=\prod_{j=0}^\infty(1-aq^j)(1-bq^j).$$
The pairing gives a meaning to $P$ as an element of $(\fr/I)^\ast$.
\end{remark}

\subsection{$6j$-symbols}

Next we turn to the problem of computing $\lx t_{kj}^M,t_{ml}^N\rx$. 
This gives the expression for our pairing in the Peter--Weyl basis, or
equivalently the action of the basis elements
 in the representation $\pi_N^+$.
Using \eqref{tx} together with  \eqref{mex} and Lemma \ref{tll} to compute 
 $\lx t_{kj}^M,t_{ml}^N\rx$, it is clear that we get a single sum, which
turns out to be a terminating ${}_8W_7$ \cite{gr}, or equivalently a
 quantum $6j$-symbol \cite{kir} or $q$-Racah polynomial \cite{aw}. 
In view of Remark \ref{rc}, we may alternatively
 deduce this  from \cite{kr}, where
the corresponding result for the representations $\pi^\om$ was derived.

Due to the large symmetry group of quantum $6j$-symbols,  there is in fact 
a large
number of single sum expressions. Omitting the details of the derivation, 
we write down one such expression  which exhibits the symmetry of the pairing.
We have checked that both methods indicated above yield the same result.
  
\begin{theorem}\label{tpt} One has the identity
$\lx t_{kj}^M,t_{ml}^N\rx=\de_{k+l,j+m=L}\,f\,T_{M+N-2L}$, where
\begin{multline*}f(\la)=
(-1)^{j+l+L}q^{(j+l-L)(\la+1+j+l-M-N)-2km}
\frac{(q^2;q^2)_{L}(q^2;q^2)_{M-k}(q^2;q^2)_{N-m}}
{(q^2;q^2)_{j}(q^2;q^2)_{l}(q^2;q^2)_{M-j}(q^2;q^2)_{N-l}}\\
\begin{split}&\times\frac{1}{1-q^{2(\la+1)}}
\frac{(q^{2(\la+1-L)};q^2)_{M-k}(q^{2(\la+1-L)};q^2)_{N-m}}
{(q^{2(\la+1-L)};q^2)_{j}(q^{2(\la+1-L)};q^2)_{l}} 
\frac{(q^{2(\la+1+M+N-2L)};q^2)_{L+1}}
{(q^{2(\la+2)};q^2)_{M-k}(q^{2(\la+2)};q^2)_{N-m}}\\
&\times\,{}_4\phi_3\left[\begin{matrix}
q^{-2k},q^{-2m},q^{2(k-M-\la-1)},q^{2(m-N-\la-1)}\\q^{-2L},
q^{2(L-M-N-\la-1)},q^{-2\la}\end{matrix};q^2,q^2\right].
\end{split}\end{multline*}
\end{theorem}

The ${}_4\phi_3$-sum appearing here is defined by \cite{gr}
$${}_{4}\phi_3\left[\begin{matrix}q^{-k},q^{-m},a,b\\
c,d,e\end{matrix};q,z\right]=\sum_{n=0}^{\min(k,m)}
\frac{(q^{-k};q)_n(q^{-m};q)_n(a;q)_n(b;q)_n}
{(q;q)_n(c;q)_n(d;q)_n(e;q)_n}\,z^n.$$

\begin{lemma}\label{pst}
One has
$$\Psi(t_{kj}^N)=C_{kj}^N(\la,\mu)\, t_{jk}^N,$$
where
$$C_{kj}^N(\la,\mu)=q^{\frac 12\la(2j-N)+\frac 12\mu(N-2k)+(j-k)(j+k+2-N)}
\frac{\qb{N}{j}}{\qb{N}{k}}\frac{(q^{2(\la+2)};q^2)_{N-j}
(q^{2(\mu+1-k)};q^2)_{N-k}}{(q^{2(\mu+2)};q^2)_{N-k}
(q^{2(\la+1-j)};q^2)_{N-j}}.$$
\end{lemma} 

This can be proved similarly as Proposition 3.12 of \cite{kr}. 

Combining
 Proposition \ref{cbp} with  Theorem \ref{tpt} and Lemma \ref{pst}
gives that the matrix 
\begin{multline*}\begin{split}& R_{lm}^{jk}(\la;M,N)=
\{t_{km}^N,t_{jl}^M\}1(\la)=\frac{1}{C_{mk}^N(\la+M-2j,\la)}\,\lx t_{mk}^N,
t_{jl}^M\rx 1(\la)\\
&\quad=\de_{j+k,l+m=L}(-1)^{j+l}q^{(l-j)(2j+2l-M-L-1)+\frac 12MN-Mk-Nj}\\
&\qquad\times\frac{(q^2;q^2)_{L}(q^2;q^2)_{M-j}}
{(q^2;q^2)_{l}(q^2;q^2)_{m}(q^2;q^2)_{M-l}}
\frac{(q^{2(\la+1-j)};q^2)_{j}(q^{2(\la+2+M+N-L-j)};q^2)_{j}}
{(q^{2(\la+2+M-2j)};q^2)_{j}(q^{2(\la+1+N-L-m)};q^2)_{l}}\\
&\qquad\times  \,{}_4\phi_3\left[\begin{matrix}
q^{-2j},q^{-2m},q^{2(j-M-\la-1)},q^{2(m-N-\la-1)}\\q^{-2L},
q^{2(L-M-N-\la-1)},q^{-2\la}\end{matrix};q^2,q^2\right]
\end{split}\end{multline*}
satisfies the dynamical Yang--Baxter equation of the form 
(here $(l_1,m_1,n_1,l_2,m_2,n_2)$ corresponds to $(a,b,c,d,e,f)$ in 
\eqref{dyx})
\begin{equation*}\begin{split}&\quad
\sum_{y=\max(0,\,m_1+n_1-N,l_2+m_2-L)}^{\min(M,\,m_1+n_1,l_2+m_2)}
R_{l_2m_2}^{l_2+m_2-y,y}(\la-2n_2+N;L,M)\,
R_{l_2+m_2-y,n_2}^{l_1,m_1+n_1-y}(\la;L,N)\\
&\times R_{y,m_1+n_1-y}^{m_1n_1}(\la-2l_1+L;M,N)=
\sum_{y=\max(0,\,l_1+m_1-L,m_2+n_2-N)}^{\min(M,\,l_1+m_1,m_2+n_2)}
R_{m_2n_2}^{y,m_2+n_2-y}(\la;M,N)\\
&\times R_{l_2,m_2+n_2-y}^{l_1+m_1-y,n_1}(\la-2y+M;L,N)\,
R_{l_1+m_1-y,y}^{l_1m_1}(\la;L,M).\end{split}\end{equation*}
This is (for $\la$ discrete) the hexagon  identity for 
$6j$-symbols, first proved by  Wigner \cite{w} for $q=1$ and by 
 Kirillov and Reshetikhin \cite{kir} in general.
Thus,  dynamical quantum groups provide an alternative algebraic framework 
for studying $6j$-symbols (in \cite{kr} we derived the pentagon or
Biedenharn--Elliott relation  using this framework).


\begin{thebibliography}{99} 

\bibitem[AW]{aw} R.~Askey and J.~A.~Wilson, {\it A set of
orthogonal polynomials that generalize the Racah coefficients 
or $6$-$j$ symbols}, SIAM J.\ Math.\ Anal.\ 10 (1979), 1008--1016. 

\bibitem[B]{b} O.~Babelon, 
{\it Universal exchange algebra for Bloch waves and
Liouville theory}, 
Comm.\ Math.\ Phys.\ 139 (1991), 619--643.

\bibitem[BBB]{bbb}O.~Babelon, D.~Bernard and E.~Billey, 
{\it A quasi-Hopf algebra interpretation of quantum
$3$-$j$ and $6$-$j$ symbols and difference equations}, Phys.\ Lett.\ B 375 
(1996), 89--97.


\bibitem[EN]{en} P. Etingof and D. Nikshych, {\it Dynamical quantum groups
at roots of $1$}, Duke Math. J. 108 (2001), 135-168.

\bibitem[EV]{ev}P.\ Etingof and A.\ Varchenko, {\it Solutions of the quantum 
dynamical Yang--Baxter equation and dynamical quantum groups}, 
Comm.\ Math.\ Phys.\ 196 (1998), 591--640.  

 \bibitem[EV2]{ev2}P.\ Etingof and A.\ Varchenko, {\it Exchange dynamical 
quantum  groups}, Comm.\ Math.\ Phys.\
205 (1999), 19--52. 




\bibitem[F]{f} G.\ Felder, {\it Elliptic quantum groups}, XIth International
 Congress of
Mathematical Physics (Paris, 1994), 211--218, Internat.\ Press, Cambridge, 
MA, 1995. 

\bibitem[FV]{fv} G.\ Felder and A.\ Varchenko, {\it On representations of the
elliptic quantum group $E_{\tau,\eta}(\mathrm{sl}_2)$}, 
Comm.\ Math.\ Phys.\ 181
(1996), 741--761.

\bibitem[GR]{gr} G.\ Gasper and M.\ Rahman, {\it Basic Hypergeometric Series},
Cambridge University Press, Cambridge, 1990.

\bibitem[GN]{gn} J.-L.~Gervais and A.~Neveu, {\it Novel triangle relation and
absence of tachyons in Liouville string field theory}, Nucl.\ Phys.\ B 238
 (1984), 125--141.



\bibitem[JKOS]{jkos} M.\ Jimbo, S.\ Odake, H.\ Konno and 
J.\ Shiraishi, {\it Quasi-Hopf twistors for elliptic quantum groups},
Transform.\ Groups 4 (1999), 303--327.

\bibitem[K]{k}  C.\ Kassel, {\it Quantum Groups}, Springer-Verlag, New York, 
1995.


\bibitem[KiR]{kir}  A.\ N.\ Kirillov and N.\ Yu.\ Reshetikhin, 
{\it Representations of 
the algebra ${U}\sb q({\rm
sl}(2))$, $q$-orthogonal polynomials and invariants of links}, 
Infinite-dimensional Lie Algebras and Groups, 285--339, 
 World Sci. Publishing, Teaneck, NJ, 1989.



\bibitem[KR]{kr}E.\ Koelink and H.\ Rosengren, {\it Harmonic analysis
on the $\mathrm{SU}(2)$ dynamical quantum group}, Acta Appl.\ Math.,
to appear.

\bibitem[M]{m}
S.\ Majid,
{\it Braided matrix structure of the Sklyanin algebra and of the quantum 
Lorentz group}, Comm.\ Math.\ Phys.\ 156 (1993),  607--638. 

\bibitem[R]{r} H.\ Rosengren,  {\it A new quantum algebraic interpretation 
of the Askey--Wilson polynomials},
Contemp.\ Math.\ 254 (2000), 371-394.

\bibitem[W]{w} E.\ P.\ Wigner, {\it On the matrices which reduce the Kronecker
products of representations of S.\ R.\ groups} (1940), in 
L.~C.~Biedenharn and H.~Van Dam (eds.), Quantum
Theory of Angular Momentum, 87--133, Academic Press, New York, 1965.

\bibitem[X]{x} P.\ Xu, {\it Quantum groupoids},
Comm.\ Math.\ Phys.\ 216 (2001),  539--581. 

\vskip 3mm
\end{thebibliography}
\end{document}